\documentclass[aop,MSNbibl,citesort,dvips]{arximspdf}
\usepackage{mathbh}

\doi{10.1214/11-AOP738} 
\volume{41}
\issue{1}
\pubyear{2013}
\firstpage{134}
\lastpage{169}

\makeatletter

\newcommand{\bb}{\mathbb}
\renewcommand{\P}{{\mathbb P}}
\newcommand{\E}{{\mathbb E}}

\newtheorem{Lemma}{Lemma}[section]
\newtheorem{theorem}[Lemma]{Theorem}
\newtheorem{Corollary}[Lemma]{Corollary}

\newproclaim{Assumption}[Lemma]{Assumption}
\newproclaim{remark}[Lemma]{Remark}
\newproclaim{Condition}[Lemma]{Condition}
\newproclaim{example}[Lemma]{Example}

\newcommand{\definedas}{\stackrel{\Delta}{=}}
\newcommand{\equalinlaw}{\stackrel{\cal L}{=}}
\newcommand{\cal}{\mathcal}
\newcommand{\bg}{{g}}
\newcommand{\ambient}{{\widetilde M}}
\newcommand{\br}{{r}}

\newcommand{\bs}{{s}}
\newcommand{\bt}{t}
\newcommand{\be}{{e}}
\newcommand{\BH}{{H}}
\newcommand{\BX}{{X}}

\newcommand{\indic}{{\mathbh{1}}}
\newcommand{\real}{{\bb{R}}}
\newcommand{\reals}{{\mathbb R}}
\newcommand{\bbr}{\reals}
\newcommand{\vep}{\varepsilon}

\newcommand{\lt}{[\![}
\newcommand{\rt}{]\!]}

\makeatother

\begin{document}
\begin{frontmatter}

\title{High level excursion set geometry for non-Gaussian infinitely divisible random fields}
\runtitle{Excursion set geometry}

\begin{aug}
\author[A]{\fnms{Robert J.} \snm{Adler}\corref{}\thanksref{t1,t4}\ead[label=e1]{robert@ee.technion.ac.il}\ead[label=u1,url]{webee.technion.ac.il/people/adler}},
\author[B]{\fnms{Gennady} \snm{Samorodnitsky}\thanksref{t1,t2}\ead[label=e2]{gennady@orie.cornell.edu}\ead[label=u2,url]{people.orie.cornell.edu/\textasciitilde gennady/}}\\
\and
\author[C]{\fnms{Jonathan E.} \snm{Taylor}\thanksref{t1,t3}\ead[label=e3]{jonathan.taylor@stanford.edu}\ead[label=u3,url]{www-stat.stanford.edu/\textasciitilde jtaylo/}}
\runauthor{R. J. Adler, G. Samorodnitsky and J. E. Taylor}
\affiliation{Technion, Cornell University and Stanford University}
\address[A]{R. J. Adler \\
Electrical Engineering\\
Technion, Haifa 32000\\
Israel \\
\printead{e1}\\
\printead{u1}}
\address[B]{G. Samorodnitsky\\
ORIE\\
Cornell University\\
Ithaca, New York 14853\\
USA\\
\printead{e2}\\
\printead{u2}}
\address[C]{J. E. Taylor\\
Department of Statistics\\
Stanford University\\
Stanford, California 94305-4065\\
USA\\
\printead{e3}\\
\printead{u3}} 
\end{aug}

\thankstext{t1}{Supported in part by US-Israel Binational
Science Foundation, 2008262.}

\thankstext{t4}{Supported in part by AFOSR FA8655-11-1-3039.}

\thankstext{t2}{Supported in part by NSA Grant MSPF-05G-049 and ARO
Grant W911NF-07-1-0078 at Cornell University.}

\thankstext{t3}{Supported in part by NSF Grant DMS-04-05970,
and the Natural Sciences and Engineering Research Council of Canada.}

\received{\smonth{7} \syear{2009}}
\revised{\smonth{9} \syear{2011}}

%
\begin{abstract}
We consider smooth, infinitely divisible random fields
$( X(\bt), \bt\in M)$, $M\subset\real^d$, with
regularly varying
L\'evy measure, and are interested in the
geometric characteristics of the excursion sets
\[
A_u = \{ \bt\in M\dvtx X(\bt) >u\}
\]
over high levels $u$.

For a large class of such random fields, we compute the $u\to\infty$
asymptotic
joint distribution of the numbers of critical points, of various types,
of $X$ in $A_u$, conditional on $A_u$ being nonempty.
This allows us, for example, to obtain the
asymptotic conditional distribution of the Euler characteristic of the
excursion set.

In a significant departure from the Gaussian situation,
the high level excursion sets for these random fields can have quite a
complicated geometry. Whereas in the Gaussian
case nonempty excursion sets are, with high probability, roughly
ellipsoidal, in the more general infinitely divisible setting almost
any shape is possible.
\end{abstract}

%

\end{frontmatter}

\section{Introduction}\label{secintro}

Let $( X(\bt), \bt\in M)$, where $M$ is a compact set
in $\real^d$
of a kind to be specified later,
be a smooth infinitely divisible random
field. We shall assume, again in a sense that we shall
make precise later, that $X$ has regularly varying tails. Note that
this means that the
tails of $X$ are heavier than exponential and, in
particular, heavier than those of a Gaussian random
field. Nevertheless,
the model we are considering allows both heavy tails (e.g., infinite
mean or variance) and light tails, in the sense of
the existence of finite moments of arbitrary given order.

We are interested in studying the excursions of the random field
over levels $u>0$, particularly when the level $u$ becomes high. Writing
%
%
\begin{equation} \label{eexcursionset}
A_u \equiv A_u(X,M) \definedas\{ \bt\in M\dvtx X(\bt) >u\}
\end{equation}
for the excursion set of $X$ over the level $u$, we shall study the geometric
characteristics of $A_u$ under the condition that it is not empty,
that is, under the condition that the level $u$ is, in fact,
exceeded.
In particular, we shall be interested in
computing the conditional limit distribution of the
Euler characteristic of $A_u$ as
$u\to\infty$. We refer the reader to~\cite{adlertaylor2007} for a
recent detailed exposition of the geometric theory of the excursion
sets of
smooth Gaussian and related random fields, and to
\cite{adlertaylorworsley2010} for applications of the theory.

In a significant departure from the well-understood Gaussian situation,
the excursion sets over high levels for the random fields in this paper
can have quite a complicated geometry. In the Gaussian
case excursion sets, unless they are empty, tend, with high
probability, to
contain a single component which is almost ellipsoidal in shape, and so
have an Euler characteristic equal to one. In
contrast, the Euler characteristics of the excursion sets in our fields
can have highly nondegenerate conditional distributions.
As a consequence, these models are sufficiently
flexible to open the possibility of fitting empirically observed excursion
sets with widely different geometric characteristics. This, more statistical,
problem is something we plan to tackle in the future.

The main result of the paper is Theorem~\ref{tmain}. While it is
rather too
technical to summarize here in full, here is the beginning of a special case.
Suppose that $N_X(i,u)$ is the number of critical points of $X$ in $A_u$
of index $i$. Thus, if $d=2$, $N_X(0,u)$ is the number of local minima of
$X$ above the level $u$ in the interior of $M$,
$N_X(1,u)$ the number of saddle points
and $N_X(2,u)$ the number of local maxima, all above the level $u$.
Then Theorem~\ref{tmain} gives
an explicit expression for the limiting joint distribution
%
%
\begin{equation}
\label{intromain}
\lim_{u\to\infty} \P\{ N_X(i,u)=n_i, i=0,\ldots,d,
|A_u\neq\varnothing\},
\end{equation}
when $M$ is the unit cube $I_d\definedas[0,1]^d$.

In fact, Theorem~\ref{tmain} goes far beyond this, since it includes not
only these critical points, but also the critical points of $X$ restricted
to the various boundaries of $I_d$ (i.e., faces, edges, etc.). The
importance of this result lies
in the fact that Morse theory shows how to use the full collection
of these critical points to describe much of the geometry of $A_u$,
whether this geometry be algebraic, integral, or differential.

Furthermore, Theorem~\ref{tmain} can also be exploited to describe a
very simple stochastic model for high
level excursion sets, as well as to develop a simple algorithm for simulating
them.

An important point to note is that although Theorem~\ref{tmain} is
stated only for $M$ the unit cube, it is ``obvious'' from the proof
that the result holds in much higher generality. For example, only
trivial changes to the proof are needed to establish the result for
convex polytopes. A little more effort will establish a version for
convex $M$ with smooth boundary. We also claim---without proof---that
Theorem~\ref{tmain}, properly reformulated, will continue to hold for
locally convex, $C^2$, Whitney stratified manifolds of the kind treated
in~\cite{adlertaylor2007}. However, in this case the additional
details that would need to be added to provide a complete proof would
take more space than justified. Thus, while we shall continue to write
$M$ for our parameter set, indicating a level of generality, throughout
the remainder of this paper we shall treat only the case $M=I_d$.

The remainder of the paper begins in Section~\ref{secthemodel}, where
we define our model, discuss the
smoothness assumptions we are imposing, as well as those related to the
regular variation
of the tails. Section~\ref{seclimitEuler} contains the main result
of the paper, on the joint distribution of the numbers of high level
critical points
of infinitely divisible random field's. This is followed with one of
its main consequences,
the distribution
of the Euler characteristic of high level excursion sets, in
Section~\ref{seclimitEuler2}. In Section~\ref{secexamples}, we
introduce a class of moving average infinitely divisible random fields
and derive conditions
under which the main result of the Section~\ref{seclimitEuler}
applies to them. We also provide examples to show that, by choosing
appropriately the parameters of the model, one can make the geometric
structure of the high level excursion sets either ``Gaussian-like'' or
``non-Gaussian-like.''
Finally, Section
\ref{secproof} contains the proof of the main theorem.

Throughout the paper, $C$ stands for finite positive constants whose
precise value is irrelevant and which may change from line to line.

\section{Smooth infinitely divisible random fields and regular variation}
\label{secthemodel}

In this section, we shall define the random fields of interest to us, describe
their distributional structure, and then specify the smoothness assumptions
necessary for studying the geometry of their excursion sets.

A reader familiar with the theory of infinitely divisible processes
will note that the route we take goes back to first principles to some
extent (e.g., it would be more standard, nowadays, to start with the
function space L\'evy measure $\lambda_X$ of
Section~\ref{lambdasection} rather than invest a couple of pages in
defining it). The need for this, as should become clear below, is to be
able to carefully define random fields, along with their first and
second order partial derivatives, on a common probability space.

\subsection{Probabilistic structure of infinitely divisible random fields}
As a first step, we shall need to define our random fields on a region
slightly larger than the basic parameter space $M$, and so, in a
notation that
will remain fixed throughout the paper, we take
$\ambient$ be a bounded open set in $\bbr^{d}$, with $M\subset
\ambient$.

We now consider infinitely divisible random fields of the form
%
%
\begin{equation} \label{eidfield}
X(\bt) = \int_S f(s;\bt) \mu(ds),\qquad t\in\ambient,
\end{equation}
where $(S, {\mathcal S})$ is a measurable space and $\mu$ is an
infinitely divisible random measure
on $S$ with characteristics defined below. (We refer you to
\cite{rajputrosinski1989} for more information on
infinitely divisible random measure s and stochastic integrals with
respect to these measures.) 

The infinitely divisible random measure $\mu$, which we shall define in
a moment, is characterized
by its
``generating triple'' $( \gamma, F, \beta)$. Here,
$\gamma$ is a $\sigma$-finite measure on $(S, {\mathcal S})$, and
plays the role of the variance measure for the Gaussian part of
$\mu$. More important for us is the L\'evy measure $F$, which is a
$\sigma$-finite measure on $S\times( \bbr\setminus
\{0\})$, equipped with the product $\sigma$-field. Finally,
$\beta$ is a signed measure on $(S, {\mathcal S})$, which plays the
role of the shift measure for $\mu$. Denote by ${\mathcal
S}_0$ the collection of sets $B$ in $\mathcal S$ for which
\[
\gamma(B) + \|\beta\|(B) + \int_{\bbr\setminus\{0\}} \lt x\rt^2
F(B,dx)<\infty,
\]
where $\|\beta\|$ is the total variation norm of $\beta$ and
\[
\lt x\rt= \cases{
x, &\quad if $|x|\leq1$,\cr
\operatorname{sign}(x), &\quad otherwise.}
\]

With all elements of the triple defined, we can now define the
infinitely divisible random measure
$(\mu(B), B\in{\mathcal S}_0)$ as a stochastic
process for which, for every sequence of disjoint ${\mathcal
S}_0$-sets $B_1,B_2,\ldots\,$, the random variables
$\mu(B_1)$, $\mu(B_2),\ldots$ are independent (i.e.,
$\mu$ is independently scattered) and if, in addition, $\bigcup_n
B_n\in{\mathcal
S}_0$, then $\mu( \bigcup_n B_n) =$ $ \sum_n \mu(B_n)$
a.s. (i.e., $\mu$ is $\sigma$-additive). Finally, for
every $B\in{\mathcal S}_0$, $\mu(B)$ is an infinitely divisible random
variable with
characteristic function given by
\[
\E\bigl\{e^{i\theta\mu(B)}\bigr\}
= \exp\biggl\{ -\frac12 \gamma(B)\theta^2
+ \int_{\bbr\setminus\{0\}} ( e^{i\theta x} -1 -i\theta\lt
x\rt) F(B,dx) + i\theta\beta(B)\biggr\}
\]
for $\theta\in\bbr$. The monograph~\cite{sato1999} can be
consulted for
information on infinitely divisible random variable's.

We shall assume (without loss of generality) that the L\'evy measure
$F$ 
has the form
%
%
\begin{equation} \label{elevymdecompose}
F(A) = \int_S \rho( s; A_s) m(ds),
\end{equation}
for each measurable $A\subset S\times( \bbr\setminus
\{0\})$, where $A_s=\{ x\in\bbr\setminus
\{0\}\dvtx(s,x)\in A\}$ is the $s$-section of the set $A$.
In
(\ref{elevymdecompose}), $m$ is a $\sigma$-finite measure on $(S,
{\mathcal S})$ (the \textit{control measure} of $\mu$), and
the measures $(\rho(s;\cdot))$
(the \textit{local L\'evy measures})
form a family of
L\'evy measures on $\reals$ such that for every Borel set $C\subset
\bbr\setminus\{0\}$, $s\to\rho(s;C)$ is a measurable function on
$S$. We can, and shall, choose the control measure $m$ in
(\ref{elevymdecompose}) in such a way that $\|\beta\|$ is
absolutely continuous with respect to $m$, and define the Radon--Nikodym
derivative $b=d\beta/dm$. The local L\'evy measures $\rho$,
which, intuitively, control the Poisson structure of the random
measure $\mu$ around different points of the space $S$, will play a
central role in all that follows.

Note that while it is possible, and common, to choose $m$ in
with the added feature that
$\gamma$ is also absolutely continuous with respect to $m$, and that
$\rho( s; \bbr\setminus\{0\})>0$ on a set of $s\in S$ of
full measure $m$, we shall not require this and so shall not do so.

Finally, we assume that the kernel $f(s;\bt), s\in S, \bt\in
\ambient$, in
(\ref{eidfield}) is deterministic and real, such that, for every
$\bt\in\ambient$, the mapping $f(\cdot; \bt)\dvtx S\to\bbr$ is
measurable, and that the
following three inequalities hold:
%
%
\begin{eqnarray} \label{eintegrgauss}
\int_S f(s;\bt)^2 \gamma(ds) &<&\infty,\\
\label{eintegrpoiss}
\int_S \int_{\bbr\setminus\{0\}} \lt xf(s;\bt)\rt^2 F(ds,dx)
&<& \infty
\end{eqnarray}
and
%
%
\begin{equation} \label{eintegrshift}\qquad
\int_S\biggl| b(s)f(s;\bt) + \int_{\bbr\setminus\{0\}} \bigl(
\lt
xf(s;\bt)\rt- \lt x\rt f(s;\bt)\bigr) \rho(s;dx)\biggr|
m(ds)<\infty.
\end{equation}
These conditions guarantee that the random field $(X(\bt), \bt
\in
\ambient)$ in (\ref{eidfield}) is well defined.

A particularly simple, but rather useful, example of this setup is
studied in Section~\ref{secexamples} below, when $X$ is a moving average
random field. In this example, both $\gamma$ and $\beta$ components of
the generating triple vanish, so, in particular, the random field
has no Gaussian component. Furthermore, $S=\reals^{d}$, the control
measure $m$ is Lebesgue,
and the local L\'evy measures $\rho(\bs,\cdot)$ are independent of
$\bs$.
Finally, the kernel function $f$ is of the form $f(\bs,\bt)=g(\bs
+\bt)$ for
some suitable $g$, and so the random field is given by
%
%
\begin{equation}
\label{2movingaverageequn}
X(\bt) = \int_{\reals^{d}} g(\bs+\bt) \mu(d\bs),\qquad\bt\in
\ambient\subset\reals^{d}.
\end{equation}
The random measure
$\mu$ has, in this case, the stationarity property
$\mu(A)\equalinlaw\mu(\bt+A)$ for all Borel $A$ of a finite Lebesgue
measure and $\bt\in\reals^{d}$, which immediately implies that a
moving average random field is stationary.
An impatient reader, who already wants to see results without wading
through technicalities, might want to now skip directly to Section
\ref{kernelssubsection}
to see what our results have to say for moving averages.

Returning to the model (\ref{eidfield}), note that it has been
defined in considerable generality, so as to
allow for as wide a range of applications as possible. For example,
we retain the Gaussian component of the random field $X$. However,
the tail assumptions imposed below will have the effect of ensuring that
the Gaussian component will not play a role in the geometric structure of
high level excursion sets.

\subsection{Regularity properties} We shall require that the sample paths
of $X$ satisfy a number of regularity properties for the theory we are
developing to hold. The main assumption will be that the paths of $X$ are
a.s. $C^2$, for which good sufficient conditions exist. The secondary
assumptions require a little more regularity, centered around the notion
of \textit{Morse functions}. For more details, including for the case of
stratified manifolds, see Chapter 9 in~\cite{adlertaylor2007}.

We need a little notation. With $M=I_d$, we write $\partial_kM$ for the
collection of the $2^{d-k}{d \choose k}$ $k$-dimensional open faces of
$M$. Thus, for example, $\partial_dM$ is the interior of~$M$, and
$\partial_0M$ the collection of $2^d$ vertices.

Next, recall that if $\ambient$ is an open neighborhood of $M$, a function
$f\dvtx\ambient\to\reals$ is called a Morse function on
$M$ if it
satisfies the following two conditions on each $\partial_k M$,
$k=0,\ldots,d$:
\begin{longlist}[(ii)]
\item[(i)]
$f_{|\partial_k M}$ it is nondegenerate on $\partial_k M$, in
the sense that the determinant of the Hessian of $f_{|\partial_k M}$
at its critical points does not vanish.
\item[(ii)] The restriction of $f$ to
$\overline{\partial_k M}=\bigsqcup_{j=0}^{k} \partial_j M$
has no critical points on
$\bigsqcup_{j=0}^{k-1} \partial_j M$.
\end{longlist}
%

Here is our first, and henceforth ubiquitous, assumption.
%
%
\begin{Assumption} \label{assc2}
On an event of probability 1, the random field $ X$ has $C^2$ sample
paths on $\ambient$
and is a Morse function on $M$.
\end{Assumption}

Sufficient conditions for Assumption~\ref{assc2} to hold are not hard
to come by. As far as the $C^2$ assumption is concerned, it suffices
to treat the Gaussian and non-Gaussian components of $X$
separately. For the Gaussian part, there is a rich and easy to apply
theory, and Section 1.4.2 of~\cite{adlertaylor2007} covers what is
needed here.

Necessary and sufficient conditions for the $C^2$ assumption on the
non-Gaussian component are not known, but a number of sufficient
conditions exist. It is not our goal in this paper to develop the best
possible conditions of this sort, so we restrict ourselves to one
situation that covers, nonetheless, a wide range of random
fields. Specifically, we shall assume that the $\gamma$ and
$\beta$ components in the generating triple of the infinitely divisible
random measure $M$ vanish,
and that the local L\'evy measures $\rho$ in
(\ref{elevymdecompose}) are symmetric; that is,
$\rho(s;-A)=\rho(s;A)$ for each $s\in S$ and each Borel
$A\in\bbr\setminus\{ 0\}$. That is, $\mu$ is a symmetric infinitely
divisible random measure without
a Gaussian component.

The following result gives
sufficient conditions for a symmetric infinitely divisible random field
without a
Gaussian component to have
sample functions in~$C^2$. The proof is not difficult, and so is left
to the reader.
(The conditions are also necessary after a slight tightening of
the assumptions on the null sets involved, cf. Theorem 5.1 of \cite
{bravermansamorodnitsky1998a}.)
%
%
\begin{theorem} \label{tc2id}
For a symmetric random field of the form (\ref{eidfield}), with $\mu$
an infinitely divisible random measure without a Gaussian
component, suppose that the kernel $f\dvtx S\times\ambient\to\bbr$ is
(product)-measurable. Assume that
for every $s\in S$ outside of set of zero $m$-measure
the function $f(s;\cdot)\dvtx\ambient\to\bbr$ is $C^2$. Furthermore,
assume that the partial derivatives
\begin{eqnarray*}
f_i(s;\bt) &=& \frac{\partial f}{\partial t_i}(s;\bt),\qquad
i=1,\ldots
, d,\\
f_{ij}(s;\bt) &=& \frac{\partial^2 f}{\partial t_i\,\partial
t_j}(s;\bt
),\qquad
i,j=1,\ldots, d,
\end{eqnarray*}
satisfy the following conditions:
\begin{longlist}[(ii)]
\item[(i)]
The integrability
condition (\ref{eintegrpoiss}) holds when the kernel $f(s;\bt)$ there
is replaced by any of the $f_i(s; \bt)$ or $f_{ij}(s; \bt)$.
\item[(ii)]
The random fields
%
%
\begin{equation}
\label{secondorders}
X_{ij}(\bt) = \int_S f_{ij}(s;\bt) \mu(ds),\qquad\bt\in
\ambient,
\end{equation}
$i,j=1,\ldots, d$, are all sample continuous.
\end{longlist}
Then the random field
$( X(t), t\in\ambient)$ has (a version with) sample
functions in $C^2$.
\end{theorem}

%
%
Thus, in searching for sufficient conditions for the a.s. second order
differentiability of $X$, it suffices to establish the continuity of the
random fields of (\ref{secondorders}).
While there are no known necessary and sufficient conditions for sample
continuity of general infinitely divisible random field's, various
sufficient conditions are
available. See, for example, Chapter 10 of
\cite{samorodnitskytaqqu1994} for the special case of stable random
fields, or~\cite{marcusrosinski2005} for some other classes of
infinitely divisible random fields.

This is as far as we shall go at the moment discussing
the issue of differentiability in Assumption~\ref{assc2}.
Conditions sufficient for $X$ to be a Morse function, also required in this
assumption, are, in principle, available as well. For example,
it follows from the arguments of Section 11.3 of~\cite{adlertaylor2007}
(cf. Theorem 11.3.1 there) that a $C^2$ field $X$ will also be, a.s.,
a Morse
function on the unit cube $I_d$
if the following two conditions are satisfied, for each face $J$ of
$I_d$, and for all $t\in J$:
\begin{longlist}[(ii)]
\item[(i)] The marginal
densities $p_t(x)$ of $\nabla X_{|J}(t)$ are continuous at 0, uniformly
in~$t$.

\item[(ii)] The conditional densities $p_t(z|x)$ of
$Z=\det\nabla^2 X_{|J}(t)$
given $\nabla X_{|J}(t)=x$ are continuous in $(z,x)$ in a neighbourhood
of 0,
uniformly in $t$.
\end{longlist}


It does not seem to be trivial to translate the above conditions into
general conditions on the kernel $f$ and the triple $(\gamma,F,\beta)$,
and we shall not attempt to do so in this paper. On the other hand,
given a specific kernel and triple, they are generally not too hard
to check. In the purely Gaussian case, simple sufficient conditions
are provided by Corollary 11.3.2 of~\cite{adlertaylor2007}, but it
is the more involved infinitely divisible case that is at the
heart of the current paper. If the latter random field is, actually a
so-called
type-$G$ random field (see~\cite{rosinski1991}) (symmetric
$\alpha$-stable random fields, $0<\alpha<2$ are a special case of
type-$G$ random fields), then these fields can be represented as
mixtures of centered Gaussian random fields, and Corollary 11.3.2 in
\cite{adlertaylor2007} may be helpful once again.



We close this section with a remark and a further assumption.
%
%
\begin{remark} \label{rkctssec}
Unless $X$ is Gaussian, Assumption~\ref{assc2} implies
that it is possible to modify the kernel $f$ in (\ref{eidfield}),
without changing the finite-dimensional distributions of $X$, in such a
way that $f(s,\cdot)$ is
$C^2$ for every $s\in
S$; see Theorem~4 of~\cite{rosinski1989}. For simplicity,
we shall therefore assume throughout that $f$
has such $C^2$ sections. This ensures, in particular, measurability of
functions of the type $\sup_{\bt\in M} |f(s,\bt)|, s\in S$,
which we shall take as given in what follows.
\end{remark}
%
%
\begin{Assumption} \label{assregkernel}
The kernel $f(s,\bt), s\in S, \bt\in\ambient$, along with its
first and
second order spatial partial derivatives $f_i$ and $f_{ij}$
are (uniformly)
bounded and, for for every $s\in S$, the
function $f(s,\cdot)$ is a Morse function on $M$.
\end{Assumption}

\subsection{The function space L\'evy measure}
\label{lambdasection} Although the infinitely divisible random field's
we are studying in this paper were constructed above via stochastic
integrals (\ref{eidfield}) and, as such, are characterised by the
triple $(\gamma,F,\beta)$ of the random measure $\mu$ and the kernel
$f$, in what follows the most important characteristic of the
infinitely divisible random field (\ref{eidfield}) will be
its\vspace*{1pt} function space L\'evy measure. This is a measure on
the cylinder sets of $\bbr^\ambient$, related to the parameters in the
integral representation of the field by the formula
%
%
\begin{equation} \label{elevyX}
\lambda_X = F\circ T_f^{-1},
\end{equation}
where $F$ is the L\'evy measure of the infinitely divisible random
measure $\mu$ and
$T_f\dvtx S\times( \bbr\setminus\{0\}) \to
\bbr^\ambient$ is given by
%
%
\begin{equation}
\label{eTf}
T_f(s,x) = xf(s,\cdot),\qquad s\in S, x\in\bbr\setminus\{0\},
\end{equation}
cf.~\cite{rajputrosinski1989}. Thus, the finite-dimensional
distributions of $X$
are given via the joint characteristic function
%
%
\begin{eqnarray} \label{ejontchf}
&&
\E\Biggl\{\exp\Biggl\{ i\sum_{j=1}^k \gamma_jX(\bt
_j)\Biggr\}\Biggr\}\nonumber\\
&&\qquad= \exp\Biggl\{ -Q(\gamma_1,\ldots, \gamma_k)
\nonumber\\[-8pt]\\[-8pt]
&&\hspace*{21pt}\qquad\quad{} + \int_{\bbr^\ambient}
\Biggl[ \exp\Biggl( i\sum_{j=1}^k \gamma_jx(\bt_j)\Biggr) -1
-i\sum_{j=1}^k \gamma_j \lt x(\bt_j)\rt\Biggr] \lambda
_X(dx)\nonumber\\
&&\hspace*{237pt}{} +
iL(\gamma_1,\ldots, \gamma_k)\Biggr\} \nonumber
\end{eqnarray}
for $k\geq1$, $\bt_1,\ldots, \bt_k\in\ambient$, and real numbers
$\gamma_1,\ldots, \gamma_k$,
where $Q$ is a quadratic function (corresponding to the Gaussian
part of $X$), and $L$ is a linear function (corresponding to the
shift). Their exact forms are not important for us at the moment.

Note that the L\'evy measures of the first
and second order partial derivatives $ X_i$ and $ X_{ij}$ are similarly
(cf. Theorem 5.1,
\cite{bravermansamorodnitsky1998a}) given by
%
%
\begin{equation} \label{elevyXpartial}
\lambda_{X_i} = F\circ T_{f_i}^{-1},\qquad\lambda_{X_{ij}} = F\circ
T_{f_{ij}}^{-1},\qquad i,j=1,\ldots, d.
\end{equation}

\subsection{Regular variation}
We now turn to the final set of technical assumptions on our infinitely
divisible random field's,
these being related to the regular variation of their L\'evy measures,
and which we formulate in terms of the local
L\'evy measures of (\ref{elevymdecompose}). These are our
final set of assumptions, and our main results hinge on them.

Recall that a function $f$ is regularly varying at infinity, with
exponent $\alpha$, if
%
%
\begin{equation}
\label{regular-variation}
\lim_{x\to\infty} \frac{f(\lambda x)}{f(x)} = \lambda^\alpha
\qquad\mbox{for all $\lambda>0$}.
\end{equation}

%
\begin{Assumption} \label{assregvar}
There exists a $H\dvtx(0,\infty)\to(0,\infty)$ that is
regularly varying at infinity with exponent $-\alpha, \alpha>
0$, and
nonnegative measurable functions $w_+$ and $w_-$ on $S$ such that
%
%
\begin{equation} \label{eregvarrho}
\lim_{u\to\infty} \frac{\rho( s;(u,\infty))}{H(u)} =
w_+(s),\qquad
\lim_{u\to\infty} \frac{\rho( s;(-\infty, -u)
)}{H(u)} =
w_-(s)
\end{equation}
for all $s\in S$. Furthermore, the convergence is uniform in the
sense there is $u_0>0$ such that, for all $u>u_0$ and all $s\in S$,
\[
\frac{\rho( s;(u,\infty))}{H(u)} \leq2 w_+(s),\qquad
\frac{\rho( s;(-\infty, -u))}{H(u)} \leq2w_-(s).
\]
\end{Assumption}

The following simple lemma relates Assumption~\ref{assregvar} to the
corresponding behaviour of the L\'evy measure $\lambda_X$ on a set of
crucial importance to us. We adopt the standard notation $a_+=\max(a,0)$
and $a_-=(-a)_+$ for the positive and negative parts of a real.
%
%
\begin{Lemma} \label{lreglevy}
Let Assumption~\ref{assregvar} hold:

\begin{longlist}[(iii)]
\item[(i)]
Assume that the kernel
$f(s,\bt), \bt\in\ambient$ is uniformly (in
$s\in S$) bounded, and that for some $\epsilon>0$,
%
%
\begin{equation} \label{esmallmaxf}
\int_S \bigl( w_+(s)+ w_-(s)\bigr) \sup_{\bt\in M}
|f(s,\bt)|^{\alpha-\epsilon} m(ds)<\infty.
\end{equation}
Then
%
%
\begin{eqnarray}
\label{ereglevy}
&&
\lim_{u\to\infty} \frac{\P\{\sup_{\bt\in M} X(\bt)
>u\}}{H(u)}\nonumber\\
&&\qquad=
\lim_{u\to\infty}\frac{\lambda_X\{ g\dvtx\sup_{\bt\in M}
g(\bt)>u\}}{H(u)}
\\
&&\qquad= \int_S \Bigl[ w_+(s)\sup_{\bt\in M}
f(s,\bt)_+^{\alpha} + w_-(s)\sup_{\bt\in M}
f(s,\bt)_-^{\alpha}\Bigr] m(ds), \nonumber
\end{eqnarray}
where $M$ can be replaced with $\ambient$ throughout. Furthermore,
%
%
\begin{eqnarray} \label{ereglevy1}\quad
\lim_{u\to\infty} \frac{\P\{\sup_{\bt\in\ambient} |X(\bt
)| >u\}}{H(u)}
&=&\lim_{u\to\infty}\frac{\lambda_X\{ g\dvtx\sup_{\bt\in
\ambient}
|g(\bt)|>u\}}{H(u)} \nonumber\\[-8pt]\\[-8pt]
&=& \int
_S \bigl( w_+(s)+ w_-(s)\bigr) \sup_{\bt\in\ambient}
|f(s,\bt)|^{\alpha} m(ds).\nonumber
\end{eqnarray}

\item[(ii)] Assume that the first order partial derivatives
$f_i(s,\bt), \bt\in\ambient$, $i=1,\ldots, d$, are uniformly (in
$s\in S$) bounded, and that for some $\epsilon>0$,
%
%
\begin{equation} \label{esmallmaxf1}
\int_S \bigl( w_+(s)+ w_-(s)\bigr) \sup_{\bt\in\ambient}
|f_i(s,\bt)|^{\alpha-\epsilon} m(ds)<\infty.
\end{equation}
Then
%
%
\begin{eqnarray} \label{ereglevyp1}
&&\lim_{u\to\infty}\frac{\lambda_{X_i}\{ g\dvtx\sup_{\bt\in
\ambient}
|g(\bt)|>u\}}{H(u)}
\nonumber\\[-8pt]\\[-8pt]
&&\qquad= \int_S \bigl( w_+(s)+ w_-(s)\bigr) \sup_{\bt\in
\ambient}
|f_i(s,\bt)|^{\alpha} m(ds). \nonumber
\end{eqnarray}

\item[(iii)] Assume that the second order partial derivatives
$f_{ij}(s,\bt), \bt\in\ambient$, $i,j=1,\ldots, d$, are uniformly (in
$s\in S$) bounded, and that for some $\epsilon>0$,
%
%
\begin{equation} \label{esmallmaxf2}
\int_S \bigl( w_+(s)+ w_-(s)\bigr) \sup_{\bt\in\ambient}
|f_{ij}(s,\bt)|^{\alpha-\epsilon} m(ds)<\infty.
\end{equation}
Then
%
%
\begin{eqnarray} \label{ereglevyp2}
&&\lim_{u\to\infty}\frac{\lambda_{X_{ij}}\{ g\dvtx\sup_{\bt
\in\ambient}
|g(\bt)|>u\}}{H(u)}
\nonumber\\[-8pt]\\[-8pt]
&&\qquad= \int_S \bigl( w_+(s)+ w_-(s)\bigr) \sup_{\bt\in
\ambient}
|f_{ij}(s,\bt)|^{\alpha} m(ds). \nonumber
\end{eqnarray}
\end{longlist}
\end{Lemma}
\begin{pf}
The first equality in (\ref{ereglevy}) follows from the second
equality there by Theorem 2.1 in~\cite{rosinskisamorodnitsky1993}.
As for the second equality in (\ref{ereglevy}), it follows from
(\ref{elevyX}) and (\ref{elevymdecompose}) that
\begin{eqnarray*}
\lambda_X\Bigl\{ g\dvtx\sup_{\bt\in\ambient}
g(\bt)>u\Bigr\}
&=&
\int_S \biggl[ \rho\biggl( s; \biggl( \frac{u}{\sup_{\bt\in
\ambient}
f(s,\bt)_+},\infty\biggr)\biggr)
\\
&&\hspace*{14.6pt}{} +
\rho\biggl( s; \biggl( -\infty, \frac{-u}{\sup_{\bt\in\ambient}
f(s,\bt)_-}\biggr)\biggr)\biggr] m(ds).
\end{eqnarray*}
Using the uniform boundedness of the kernel and Potter's bounds
(cf.~\cite{Potter} or~\cite{binghamgoldieteugels1987},
Theorem 1.5.6) 
we see that for any $\epsilon>0$ there is $C>0$
such that for all $u>1$,
\[
\frac{ \rho( s; ( {u}/{\sup_{\bt\in\ambient}
f(s,\bt)_+},\infty))}{H(u)}
\leq C \sup_{\bt\in\ambient} f(s,\bt)_+^{\alpha-\epsilon}
\]
and
\[
\frac{ \rho( s; ( -\infty, -{u}/{\sup_{\bt\in
\ambient}
f(s,\bt)_-}))}{H(u)}
\leq C \sup_{\bt\in\ambient} f(s,\bt)_-^{\alpha-\epsilon}.
\]
The limit (\ref{ereglevy}) now follows from Assumption
\ref{assregvar} via (\ref{esmallmaxf}), regular variation,
and dominated convergence.
The proof of (\ref{ereglevy1}) is identical, as are the
proofs of (ii) and~(iii).
\end{pf}
%
%
\begin{remark}
The assumption of uniform boundedness of the kernel $f$ in
(\ref{eidfield}) and its partial derivatives
will be kept throughout the paper (it is already a
part of Assumption~\ref{assregkernel}), but the only
place it is used is in Lemma~\ref{lreglevy}. It is not difficult to
see that this assumption can be removed at the expense of appropriate
assumptions on the behaviour near the origin of the local
L\'evy measures in (\ref{elevymdecompose}) and of slightly
modifying the integrability condition
(\ref{esmallmaxf}). Given that this paper is already rather
heavy on notation, we shall continue to work with uniform
integrability, which helps keep things comparatively tidy.
Note that it is also clear that, for the purpose of proving
(\ref{ereglevy}) alone, the integrability assumption
(\ref{esmallmaxf}) could be relaxed.
\end{remark}

\section{Limiting distributions for critical points}
\label{seclimitEuler}

Our initial aim, as described in the \hyperref
[secintro]{Introduction}, was
to obtain information about the distribution of the Euler
characteristic of the excursion sets of (\ref{eexcursionset}). As is
known from Morse critical point theory, Euler characteristics of
excursion sets are closely related to the critical points above fixed
levels. We shall describe this connection in the following section and,
for the moment, concentrate on the critical points of $X$, which are
also of intrinsic interest.

Recall the partition of the cube $M$ into collections $\partial_k M$
of facets of dimension~$k$. Let $J$
denote one such facet, of dimension $0\leq k\leq d$.

Let $g$ be a $C^2$ function on $\ambient$, and
and for $i=0,1,\ldots, \dim(J)$, let
${\mathcal C}_g (J; i)$ be the set of \textit{critical points of
index} $i$ of $g_{|J}$. These are the points for which $\nabla g(t)$ is
normal to $J$ at $t$, and for which the index of
the Hessian of $g_{|J}$, computed with respect to the natural
orthonormal basis of $J$ and
when considered as a matrix, has index $i$.
(Recall that the index of a matrix is the number of its
negative eigenvalues.)
Let
\[
N_g(J;i) = \operatorname{Card} ( {\mathcal C}_g(J;i))
\]
and, for real $u$,
\[
N_g(J;i:u) = \operatorname{Card} \bigl( {\mathcal C}_g(J;i)\cap\{
\bt\dvtx
g(t)>u\}\bigr)
\]
be the the overall number of the critical points of different types of
$g_{|J}$, and
the number of these critical points above the level $u$, correspondingly.
Since $g$ is a Morse function, it is standard fare that
all of the above numbers are finite
(e.g.,~\cite{adlertaylor2007}).

Just a little more notation is required for the main theorem. Let
$f$ be the kernel in the integral representation (\ref{eidfield}) of
an infinitely divisible random field. For $k=0,1,\ldots, d$, a facet
$J$ and
$i=0,1,\ldots, \dim(J)$, let
\[
c_i(J;s) = N_{f(s;\cdot)}(J;i)
\]
be the number
of the critical points of the $s$-section of $f$ of the appropriate
type, well defined since by Assumption~\ref{assregkernel} the
sections are Morse functions.

Furthermore, let $( \bt_l(J;i;s), l=1,\ldots,
c_i(J;s))$ be an enumeration of these critical points, and,
for $1\leq m\leq c_i(J;s)$ let
\[
f_{[m]}^{(J;i:+)}(s),\qquad
f_{[m]}^{(J;i:-)}(s)
\]
be, correspondingly, the $m$th largest of the positive and negative
parts of
$f(s;\bt_l(J; i;s))$.
[Both quantities are set to zero if $m>c_i(J;s)$.]

Finally, extend
these definitions to $m=0$ by setting
\[
f_{[0]}^{(J;i:+)}(s)
= \sup_{\bt\in M}(f(s;\bt))_+,\qquad
f_{[0]}^{(J;i:-)}(s)
= \sup_{\bt\in M}(f(s;\bt))_-.
\]

The following theorem, proven in Section~\ref{secproof}, is the main
result of this paper. It describes the limiting, conditional, joint
distribution of the number of critical points of all possible types of
a infinitely divisible random field over the level $u$, as
$u\to\infty$, given that the random field actually exceeds level $u$ at
some point. We recall that, in the theorem, $M$ is the unit cube and
$\ambient$ an open, bounded, neighborhood of $M$. However, extensions
to more general $M$, as explained in the
\hyperref[secintro]{Introduction}, can also be proven.
%
%
\begin{theorem} \label{tmain}
Let $( X(\bt), \bt\in\ambient)$ be an infinitely divisible random field
with representation (\ref{eidfield}), satisfying Assumptions~\ref{assc2},
\ref{assregkernel} and~\ref{assregvar}. Assume that
(\ref{esmallmaxf}) holds for some $\epsilon>0$. Then, for any
collection $\cal J$ of facets
in the various $\partial_kM$, $k\in\{0,1,\ldots,d\}$,
and any collection of nonnegative integers
$\{n(J;i)=0,1,\ldots, i=0,1,\ldots, \dim(J), J\in\cal J\}$, as
$u\to\infty$,
%
%
\begin{eqnarray}
\label{elimdist}
&& \P\Bigl\{ N_{X}(J;i:u)\geq n(J;i), J\in{\mathcal J},
i=0,1,\ldots, \dim(J) \big| \sup_{\bt\in M}
X(\bt)>u\Bigr\} \nonumber\\
&&\qquad\to{\int_S \Bigl[ w_+(s)\Bigl( \min_{J,i}
f_{[n(J;i)]}^{(J;i:+)}(s)\Bigr)^\alpha+ w_-(s)\Bigl( \min_{J,i}
f_{[n(J;i)]}^{(J;k-i:-)}(s)\Bigr)^\alpha\Bigr] m(ds)}\\
&&\qquad\quad\hspace*{0pt}{}\times\biggl({\int_S
\Bigl[
w_+(s)\sup_{\bt\in M} f(s,\bt)_+^{\alpha} +
w_-(s)\sup_{\bt\in M} f(s,\bt)_-^{\alpha}\Bigr]
m(ds)}\biggr)^{-1}.\nonumber
\end{eqnarray}
\end{theorem}
%
%
\begin{remark} \label{rkmain}
While the structure of (\ref{elimdist}) might be rather forbidding at
first sight, its meaning is actually rather simple.
The main point of Theorem~\ref{tmain} is that, once the random field
reaches a high level, its behavior above that level is very similar to
that of the much simpler random field,
%
%
\begin{equation}
\label{zfield}
Z(\bt) = Vf(W,\bt),\qquad\bt\in\ambient,
\end{equation}
where $(V,W)\in(\bbr\setminus\{ 0\})\times S$ is a random pair, the
joint law of which is the finite restriction of the L\'evy measure $F$
to the
set
\[
\Bigl\{ (x,s)\in(\bbr\setminus\{ 0\})\times S\dvtx
\sup_{\bt\in\ambient} |xf(s;\bt)|>1\Bigr\},
\]
normalized to be a probability measure on that set.

\end{remark}
%
%
\begin{remark} \label{rkalgo}
In fact, one can go much further than in the previous remark, and interpret
the limit (\ref{elimdist}) as showing that
limiting conditional joint distribution of critical
points is a mixture distribution, that can be described as follows. Set
\[
H = \int_S \Bigl[
w_+(s)\sup_{\bt\in M} f(s,\bt)_+^{\alpha} +
w_-(s)\sup_{\bt\in M} f(s,\bt)_-^{\alpha}\Bigr] m(ds).\vadjust{\goodbreak}
\]

\begin{longlist}[(1)]
\item[(1)] Select a random point $W\in S$ with probability law
$\eta$ on $S$ where
\[
\frac{d\eta}{dm}(s) = H^{-1} \Bigl[ w_+(s)\sup_{\bt\in M}
f(s,\bt)_+^{\alpha} + w_-(s)\sup_{\bt\in M}
f(s,\bt)_-^{\alpha} \Bigr],\qquad s\in S.
\]
\item[(2)] Given $W=s$, select a random value $I\in\{-1,1\}$ with
the law
\[
P( I=1| W=s) =
\frac{w_+(s)\sup_{\bt\in M} f(s,\bt)_+^{\alpha}}
{w_+(s)\sup_{\bt\in M} f(s,\bt)_+^{\alpha} +
w_-(s)\sup_{\bt\in M} f(s,\bt)_-^{\alpha}}.
\]
\item[(3)] Let $V_\alpha$ be a random variable independent of $W$ and $I$,
with $\P\{V_\alpha\leq x\} = x^\alpha$ for $0\leq x\leq1$. Then
the numbers of critical points $(N_{X}(J;i:u), J\in
{\mathcal J})$, given that $\sup X>u$, have, as $u\to\infty$,
the same distribution as the numbers of critical points of the random field
%
%
\begin{equation}\label{eECp}
\biggl( \frac{f(W,\bt)_+}{\sup_{\br\in M}f(W,\br)_+}, \bt\in
M\biggr)
\end{equation}
above the level $V_\alpha$ if $I=1$, and the numbers of
critical points of
%
%
\begin{equation}\label{eECm}
\biggl( \frac{f(W,\bt)_-}{\sup_{\br\in M}f(W,\br)_-}, \bt\in
M\biggr)
\end{equation}
above the level $V_\alpha$ if $I=-1$.
\end{longlist}
\end{remark}
%
%
\begin{remark}
\label{extendedremark}
While Theorem~\ref{tmain} counts critical points classified by
their indices, there are also other properties of critical points
that are of topological importance. For example, in
\cite{adlertaylor2007} considerable emphasis was laid on the
so-called ``extended outward critical points.''
These are the
critical points $t\in M$ for which $\nabla f(t) \in N_t(M)$,
where $N_t(M)$ is the normal cone of $M$ at $t$.

Extended outward critical points play a major role in Morse theory, in
terms of defining the Euler characteristics of excursion sets. It will
be easy to see from the proof of Theorem~\ref{tmain} that its
the statement remains true if one replaces critical points by extended
outward critical points. This will be used in certain applications of
Theorem~\ref{tmain} below.
\end{remark}

\section{The Euler characteristic of excursion sets}
\label{seclimitEuler2}

One application of Theorem~\ref{tmain} is to the Euler
characteristic $\varphi(A_u)$ of the excursion set $A_u$ over a
high level~$u$.
We shall not define the Euler characteristic here, but rather send you to
\cite{adlertaylor2007} for details. The Euler characteristic of an
excursion set of a Morse function is equal to the alternating sum of
the numbers of extended outward critical points of the function over
the level. This leads to the following result, an immediate corollary
of Theorem~\ref{tmain}, (9.4.1) in~\cite{adlertaylor2007}, and
Remarks~\ref{rkalgo} and~\ref{extendedremark} above.
%
%
\begin{Corollary} \label{ceec}
Under the conditions of Theorem~\ref{tmain}, the conditional
distribution of the Euler characteristic of the excursion set of an
infinitely divisible
random field computed with its limiting conditional distribution given\vadjust{\goodbreak}
that the level is exceeded, is given by the mixture of the Euler
characteristics of the random fields (\ref{eECp}) and
(\ref{eECm}), with the mixing distribution as described in Remark
\ref{rkalgo}. In particular, the expected Euler
characteristic of the excursion set of the limiting (conditional)
random field is given by
\[
H^{-1} \int_S \Bigl[ w_+(s)\sup_{\bt\in M}f(s,\bt)_+^{\alpha} \E
\{C_+(s)\}
+ w_-(s)\sup_{\bt\in M} f(s,\bt)_-^{\alpha}\E\{C_-(s)\} \Bigr]
m(ds).
\]
Here, for $s\in S$, $C_\pm(s)$ is the Euler characteristic of the
excursion set of the field $( f(s,\bt)_\pm/\sup_{\br\in
M}f(s,\br)_\pm, \bt\in M)$ above the level
$V_\alpha$.
\end{Corollary}

\section{An example: Moving average fields}
\label{secexamples}

The power and variety of the results of the previous two sections
can already be seen in a relatively simple but application rich
class of random fields,
the moving average fields with kernel $g$ that were introduced at
(\ref{2movingaverageequn}).
Our basic assumptions, that will hold throughout this section, are:
\begin{longlist}[(iii)]
\item[(i)] The function $g$ is $C^2$ on $\bbr^d$ and satisfies
(\ref{eintegrpoiss}) and (\ref{eintegrshift}).
\item[(ii)] $\mu$ is an infinitely divisible random measure on $\bbr
^d$, for which the Gaussian
and shift components in the generating triple, $\gamma$ and $\beta$, vanish.
\item[(iii)] The control measure $m$ in
(\ref{elevymdecompose}) is $d$-dimensional Lebesgue measure.
\item[(iv)] The local L\'evy measures $\rho(\bs, \cdot) = \rho
(\cdot)$ are
independent of $\bs\in\bbr^d$.
\end{longlist}

By choosing different kernels $g$, we shall see that quite different
types of high level excursion sets arise, as opposed to the Gaussian
case, in which ellipsoidal sets are, with high probability, ubiquitous.


\subsection{\texorpdfstring{Checking the conditions of Theorem \protect\ref{tmain} for type $G$ moving averages}
{Checking the conditions of Theorem 3.1 for type G moving averages}}

In this subsection, we exhibit a broad family of moving average random
fields (\ref{2movingaverageequn}) 
for which we shall verify the
conditions required by the main result of Section
\ref{seclimitEuler}. These are the so-called type $G$ random
fields. We emphasize that the applicability of our main results is not
restricted to type $G$ random fields. For the latter, we can use
standard tools to check the assumptions of Theorem
\ref{tmain}, which is why they are presented here. The main result of
this subsection is the following theorem.
%
%
\begin{theorem} \label{csatisfies}
A moving average infinitely divisible random field $\BX$ satisfying
Conditions~\ref{conrho} and
\ref{cong} below also satisfies the assumptions of Theorem~\ref{tmain}.
\end{theorem}
%
%
\begin{Condition} \label{conrho}
The local L\'evy measure $\rho$ is a symmetric measure of the form
%
%
\begin{equation} \label{esubG}
\rho(B) = \E\{ \rho_0(Z^{-1}B)\},
\end{equation}
where $B$ is a Borel set,
$Z$ is a standard normal random variable, and $\rho_0$ is a
symmetric L\'evy measure on $\bbr$. Furthermore,
the function $\rho_0( (u,\infty)), u>0$
is regularly varying at infinity with exponent $-\alpha, \alpha>
1$, and\vadjust{\goodbreak} there is $\beta\in[1,2)$ such that
\[
\rho_0( (u,\infty))\leq au^{-\beta}
\]
for all $0<u< 1$, for some $0<a<\infty$.
\end{Condition}

In fact, for \textit{any} L\'evy measure $\rho_0$ on $\bbr$,
(\ref{esubG}) defines a L\'evy measure on $\bbr$; see,
for example, Proposition 2.2 in~\cite{maejimarosinski2002}.
Furthermore, it
is simple to check that the behavior of the measures $\rho$ and $\rho_0$
are similar at zero and infinity. Specifically,
%
%
\begin{equation} \label{erhoinf}
\lim_{u\to\infty} \frac{\rho( (u,\infty))}{\rho
_0(
(u,\infty))} = \E\{Z_+^\alpha\}
\end{equation}
and
%
%
\begin{equation} \label{erhozero}
\rho( (u,\infty)) \leq\E\{\max( |Z|^\beta,
1)au^{-\beta}\}
\end{equation}
for $0<u< 1$. In particular, a moving average infinitely divisible
random field satisfying
Condition
\ref{conrho} automatically also satisfies Assumption
\ref{assregvar}. It suffices to choose $H(u) = \rho_0(
(u,\infty)), u>0$, and $w_+(\bs) = w_-(\bs) = \E\{
Z_+^\alpha\}$.

It is Condition~\ref{conrho} that makes the random field a ``type $G$
random field.'' It implies that the random field $X$ can be represented
as a certain mixture of stationary Gaussian fields, cf.
\cite{maejimarosinski2002}. Under the conditions we impose, each one
of the latter is a.s. a Morse function, which will tell us that the
moving average itself has sample functions which are, with
probability 1, Morse functions.

If it is known from
other considerations that the sample functions of a specific infinitely
divisible random field
are, with
probability 1, Morse functions, then (\ref{esubG}) is not needed, and
only the assumptions on the behavior of the tails of the L\'evy measure
$\rho( (u,\infty))$ as $u\to0$ or $u\to\infty$
are required. In the present form of Condition~\ref{conrho}, these
assumptions become the conclusions (\ref{erhoinf}) and
(\ref{erhozero}) from the corresponding assumptions on the L\'evy
measure $\rho_0$.

%
\begin{Condition} \label{cong}
The kernel $g$ is in $C^3$, and its restriction to any bounded
hypercube is a Morse function. Assume that the first and the second
derivatives $g_i, i=1,\ldots, d$, and $g_{ij}, i,j=1,\ldots, d$,
satisfy (\ref{eintegrpoiss}) and (\ref{eintegrshift}). Assume,
further, that for almost every
$\bs\in\bbr^d$ there is no subspace of dimension strictly less than
$(d^2+3d)/2$ to which the vectors $( g_i(\bs), i=1,\ldots, d,
g_{ij}(\bs), i,j=1,\ldots, d$, $i\leq j)$ belong.

Finally, assume that the function
\[
T_g(\bs) = \sup_{\bt\in[-1,1]^d} |g(\bs+\bt)|,\qquad\bs\in\bbr^d,
\]
satisfies $T_g\in L^{\alpha-\vep}(\bbr^d)$, while the function
\[
\widetilde T_g(\bs) = \max_{i,j\in1,\ldots, d} |g_{ij}(\bs)|
+
\sup_{\bt\in[-1,1]^d, i,j,k\in1,\ldots, d}
|g_{ijk}(\bs+\bt)|,\qquad\bs\in\bbr^d,
\]
satisfies $\widetilde T_g\in L^{\alpha-\vep}(\bbr^d)\cap L^\beta
(\bbr
^d)$ for
some $\vep>0$ and for the $\alpha$ and $\beta$ for which Condition
\ref{conrho} holds.
\end{Condition}

Since these are assumptions on the kernel $g$ in the integral
representation (\ref{2movingaverageequn}) 
of the random field, and
the kernel is often explicitly given, the above conditions are,
generally, easy to apply. See the examples below.

Clearly, a moving average infinitely divisible random field satisfying
Condition~\ref{cong} will also satisfy Assumption~\ref{assregkernel}.
It also satisfies (\ref{esmallmaxf}). Theorem~\ref{csatisfies} is then
an immediate consequence of the following two lemmas, the first one of
which follows in a straightforward fashion from the metric entropy
condition in Remark~2.1 in~\cite{marcusrosinski2005}.
%
%
\begin{Lemma}
\label{continuitylemma}
A moving average infinitely divisible random field satisfying Condition
\ref{cong} has sample
paths in $C^2$.
\end{Lemma}
%

To complete the proof of Theorem~\ref{csatisfies},
we need to check that a moving average satisfying Conditions
\ref{conrho} and~\ref{cong} has sample functions that are, with
probability~1, Morse functions. As mentioned above,
we shall accomplish this by representing the
random field $\BX$ as a mixture of zero mean Gaussian random fields,
each one of which will have, with probability 1, sample functions that
are Morse functions.
%
%
\begin{Lemma} \label{prdeducec2}
A moving average infinitely divisible random field satisfying
Conditions~\ref{conrho} and
\ref{cong} has sample functions that are, with probability 1, Morse
functions.
\end{Lemma}
\begin{pf}
Let $\nu$ and $\widetilde\nu$ be probability measures on $\bbr$
and $\bbr^d$ absolutely continuous with respect to the L\'evy
measure $\rho_0$ and to $d$-dimensional Lebesgue measure~$\lambda_d$,
respectively.
Let
\[
\psi(x) = \frac{d\nu}{d\rho_0}(x),\qquad x\in\bbr,\qquad\varphi
(\bs) =
\frac{d\widetilde\nu}{d\lambda_d}(\bs),\qquad\bs\in\bbr^d.
\]
Then the random field $\BX$ has a representation as an infinite sum of
the form
%
%
\begin{equation} \label{eseriesexp}
X(\bt) = \sum_{k=1}^\infty Z_kV_k g(\bt+ \BH_k)\indic\bigl(
\psi(V_k)\varphi(\BH_k) \Gamma_k\leq1\bigr),\qquad\bt\in\bbr
^d,
\end{equation}
where $(Z_1,Z_2,\ldots)$ are i.i.d. standard normal random variables,
$(V_1,V_2,\ldots)$ are i.i.d. random variables with a common law
$\nu$, $(\BH_1,\BH_2,\ldots)$ are i.i.d. random vectors in $\bbr^d$
with a common law $\widetilde\nu$, and $(\Gamma_1,\Gamma_2,\ldots
)$ are the
points of a unit rate Poisson process on $(0,\infty)$. All four sequences
are independent. See
\cite{maejimarosinski2002}, Section~5, for details. Furthermore, by Lemma
\ref{continuitylemma} and Theorem~\ref{tc2id}, the first and
second order
partial derivatives of $\BX$ are also moving average
random fields, with corresponding series
representations. In particular, all that needs changing is to replace
$g$ in (\ref{eseriesexp})
by an appropriate derivative.\vadjust{\goodbreak}

We may assume, without loss of generality, that the standard Gaussian
sequence $(Z_1,Z_2,\ldots)$ is defined on a probability space $(
\Omega_1, {\mathcal F}_1, \P_1)$, and the remaining random
variables on the right-hand side of (\ref{eseriesexp})
are defined on a different probability space $(
\Omega_2, {\mathcal F}_2, \P_2)$, so that the random fields 
defined by the series are defined on the product probability space
$( \Omega, {\mathcal F}, \P) = (
\Omega_1\times\Omega_2, {\mathcal F}_1\times{\mathcal F}_2,
\P_1\times\P_2)$. Thus, for every fixed $\omega_2\in\Omega_2$,
the conditional random field $\BX( (\omega_1,\omega_2)),
\omega_1\in\Omega_1$, is a centered Gaussian random field. We now
apply to this random field Corollary 11.3.2 in
\cite{adlertaylor2007}.

Firstly, we check the condition on the incremental variance of the
second order partial derivatives
there. In obvious notation, for every $i,j=1,\ldots,d$, and
\mbox{$\bt, \bs\in M$},
\begin{eqnarray*}
&&\E_1\bigl\{\bigl( \BX_{ij}(\bt)- \BX_{ij}(\bs)\bigr)^2
\bigr\}
\\
&&\qquad= \sum_{k=1}^\infty V_k^2 \bigl( g_{ij}(\bt+ \BH_k) -
g_{ij}(\bs+
\BH_k)\bigr)^2 \indic\bigl(
\psi(V_k)\varphi(\BH_k) \Gamma_k\leq1\bigr).
\end{eqnarray*}
Bounding the H\"older constant of a function by the largest value of
its partial derivatives, as in the proof of Lemma
\ref{continuitylemma}, we obtain
\begin{eqnarray*}
&&\E_1\bigl\{\bigl( \BX_{ij}(\bt)- \BX_{ij}(\bs)\bigr)^2
\bigr\}
\\ &&\qquad
\leq d^2 \|\bt-\bs\|^2 \sum_{k=1}^\infty V_k^2 \widetilde T_g^2(\BH_k)
\indic\bigl(
\psi(V_k)\varphi(\BH_k) \Gamma_k\leq1\bigr)
\end{eqnarray*}
and, hence, the incremental variance condition will follow once we
check that the infinite sum above converges. For this, we need to
check (see~\cite{rosinski1990b}) that
\[
\int_0^\infty\E_2\bigl\{\min\bigl[ 1, V^2 \widetilde T_g(\BH)^2
\indic\bigl(
\psi(V)\varphi(\BH) x\leq1\bigr)\bigr]\bigr\} \,dx<\infty.
\]
(The random variables without a subscript represent generic
members of the appropriate sequences.) By the definition of the
derivatives $\psi$ and $\varphi$, this reduces to checking that
\[
\int_\bbr\int_{\bbr^d} \min[ 1, y^2 \widetilde T_g(\bs
)^2]
\,d\bs\,dy<\infty,
\]
which is an elementary consequence of the integrability assumptions
imposed on $\widetilde T_g$ in Condition~\ref{cong}, and of the
assumptions imposed on the L\'evy measure $\rho_0$ in Condition
\ref{conrho}.

It remains to check that the joint
distribution under $\P_1$ of the random vectors
of partial derivatives $( X_i, X_{ij})$ is
nondegenerate for $P_2$-almost every $\omega_2$. This, however,
follows from representing  
the derivatives as sums
akin to (\ref{eseriesexp}), along with the part of Condition
\ref{cong} that rules out the possibility that the derivatives of the
kernel $g$ belong to a lower dimensional subspace.\looseness=1
\end{pf}

\subsection{Examples: How the shape of the kernel can affect the geometry of excursion sets}
\label{kernelssubsection}

In the examples below, Condition~\ref{conrho} is a standing
assumption, and will not be mentioned explicitly. Our first example
is of an infinitely divisible moving average random field whose high
level excursion
sets have a similar geometric structure to those of Gaussian random fields.
%
%
\begin{example} \label{exlikegauss}
Let $g$ be a nonnegative kernel satisfying Condition~\ref{cong},
that is also rotationally invariant and radially decreasing; that is,
$g(\bt) = g_r(\|\bt\|)$ for some nonnegative, decreasing
$g_r$ on $[0,\infty)$. An example is the Gaussian kernel
$g(\bt) = \exp\{ -a \|\bt\|^2\}$, $a>0$, for which it is
trivial to check
that the restrictions in
Condition~\ref{cong} on the various partial derivatives of $g$ hold.

Corollary~\ref{ceec} tells us that the
Euler characteristic of the excursion set over a high level, given
that the level is exceeded, is asymptotically that of the field
%
%
\begin{equation} \label{eshiftedg}
\Bigl( \sup_{\br\in I_d} g(\bs+\br)\Bigr)^{-1}
g(\bs+\bt),\qquad \bt\in I_d,
\end{equation}
with a randomly chosen $\bs\in\bbr^d$ and over a random level
$V_\alpha$.

The assumption of rotational invariance and radial monotonicity on the
kernel $g$ implies that, in this case, the excursion set of the random
field is the intersection of a Euclidian ball centered at the point
$-\bs$ and the cube $I_d$. This is a convex set and, hence,
has Euler characteristic equal to 1, regardless of the point
$\bs\in\bbr^d$ or the random level $V_\alpha$.

In this case, the limiting conditional distribution of the Euler
characteristic is degenerate at the point 1. Furthermore,
the excursion set has, with high probability, a
``ball-like shape,'' as is the case for smooth Gaussian random fields.
\end{example}

In spite of the ``Gaussian-like'' conclusion in the previous example, it
is easy to modify it to make the high level excursion sets of an
infinitely divisible random field behave quite differently. Here is a
simple example.
%
%
\begin{example} \label{exunlikegauss}
We modify the kernel $g$ of the previous example by adding to
it oscillations, while preserving its
smoothness and integrability properties. For example, for fixed $\theta
\in\real^d$, take
\[
g(t) = ( 1+ \cos\langle\theta, t\rangle) e^{-a
\|t\|^2},\qquad t\in\bbr^d.
\]
Then, depending on the random choice of the point $\bs$ in
(\ref{eshiftedg}), the structure of the excursion sets in
$I_d$ could be quite varied, as it depends on the shape of
$g$ in the translated cube $I_d^{(-s)}$.
Thus, depending on the random level $V_\alpha$,
the shape of the excursion set may be quite
different from a ball-like shape. In particular, its
Euler characteristic will have a nondegenerate
distribution.\looseness=1
\end{example}

\subsection{The bottom line}
The bottom line, of course, is that
the shape of the excursion sets is determined, to a large extent, by
the shape of the kernel in the integral representation of the random
field or, alternatively, by the geometric properties of the functions
on which the L\'evy measure of the random field is supported. By
choosing appropriate parameters for the random field, one can generate
quite different distributions for the Euler and other geometric
characteristics of high level excursion sets.

Our hope is that this fact will generate greater flexibility in
applications, allowing the practitioner to choose
models with predetermined excursion set shapes.
Furthermore, the description of the limiting conditional
\textit{distribution} (and not only the expected value) of the numbers of
critical points and so the Euler
characteristic should allow one to devise better statistical tests
based on the observed excursion sets.

\section{\texorpdfstring{Proof of Theorem \protect\ref{tmain}}{Proof of Theorem 3.1}}
\label{secproof}

The proof is rather long and rather technical, although the basic idea
is not difficult.

The basic idea, which is common to many proofs involving infinitely
divisible random field's $X$, is
to write $X$ as a sum of two parts, one which tends to be large and one
which is made up of smaller perturbations. The large part, which,
distributionally, behaves as a Poisson sum of deterministic functions
with random multipliers, is comparatively simple to handle, and it is
this part that actually accounts for the limit in Theorem
\ref{tmain}. One then needs to show that the small pertubations can be
ignored in the $u\to\infty$ limit. In the argument that follows this
is somewhat more difficult
than is usually the case, since even if the small part is small in magnitude
it can, in principle, have a major effect on variables such as the number
of critical points of the sum. [Think of any smooth function
$f$
to which is added $g(t)=\epsilon\cos(\langle\theta, t\rangle
/\epsilon^2)$.
No matter how large
$\lambda$ might be, nor how small $\epsilon$ might be, the critical points
of $f+g$ are, effectively, determined by $g$, not~$f$.]

Due to the length of the ensuing proof, we shall do our best to
signpost it
as it progresses.

\begin{longlist}[(iii)]
\item[(i)]
\textit{Some notation for the parameter space and for
critical points.} As mentioned earlier, in this section we shall take
as our parameter
space the cube $I_d$. The first step is to develop notation for
describing its stratification.

Let ${\mathcal J}_k$ be the collection of the $2^{d-k}{d \choose k}$
faces of $I_d$ of dimension $k$, $k=0,\ldots,d$, and
let ${\cal J} = \bigcup_{k} {\cal J}_k$.
For each face $J\in
{\mathcal J}_k$, there is a corresponding set $\sigma(J)\subseteq\{
1,\ldots, d\}$ of cardinality $k$ and a sequence $\epsilon(J)\in
\{-1,1\}^{\sigma(J)}$ such that
\[
J = \{ \bt=(t_1,\ldots, t_d)\in I_d\dvtx t_j=\epsilon_j
\mbox{ if $j\notin\sigma(J)$
and } 0<t_j<1 \mbox{ if
$j\in\sigma(J)$}\}.
\]
Let $g$ be a $C^2$ function on an open set $\ambient$ containing
$I_d$. For $J\in{\mathcal J}_k$ and $i=0,1,\ldots, k$, let
${\mathcal C}_g (J; i)$ be the set of points $\bt\in J$ satisfying the
following two conditions:
%
%
\begin{equation} \label{ecriticalcond1}
\frac{\partial g}{\partial t_j}(\bt)=0 \qquad\mbox{for each
$j\in\sigma(J)$},
\end{equation}
%
%
%
\begin{eqnarray} \label{ecriticalcond3}
\nonumber\\[-18pt]
&&\mbox{the matrix } \biggl( \frac{\partial^2 g(\bt)}{\partial
t_m\,\partial
t_n}\biggr)_{m,n\in\sigma(J)} \mbox{ has nonzero determinant}
\\ &&\mbox{and its
index is equal to $k-i$.} \nonumber
\end{eqnarray}
%

Now define $N_g(J;i)$ and $N_g(J;i:u)$ in terms of ${\mathcal C}_g(J;i)$
as in Section~\ref{seclimitEuler}.


\item[(ii)] \textit{Splitting $X$ into large and small components.}
By Assumption~\ref{assc2}, $X$ and its first and
second order partial derivatives are a.s. bounded on $\ambient$,
and, by (\ref{elevyX}), the L\'evy measure of $X$
is concentrated on $C^2$ functions.
%
Defining
\[
S_L = \Bigl\{ g\in C^2\dvtx\max\Bigl[ \sup_{\bt\in\ambient}
|g(\bt)|,
\sup_{\bt\in\ambient, i=1,\ldots,d} |g_i(\bt)|,
\sup_{\bt\in\ambient, i,j=1,\ldots,d} |g_{ij}(\bt)|\Bigr]
> 1\Bigr\},
\]
the sample boundedness of $X$, along with
(\ref{elevyXpartial}) and general properties of
L\'evy measures on Banach spaces (e.g.,~\cite{linde1986})
imply that
%
%
%
\begin{equation}
\label{efinitelevy}
\theta\definedas\lambda_X\{ S_L \} < \infty.
\end{equation}

We are now ready to decompose the infinitely divisible
random field $X$ into a sum of two
independent infinitely divisible components by writing
%
%
\begin{equation}\label{erfdecompose}
X(\bt) = X^L(\bt) +Y(\bt),\qquad \bt\in\ambient,
\end{equation}
where $X^L$ is a compound Poisson random
field with characteristic functions, which,
for $k\geq1$, $\bt_1,\ldots, \bt_k\in\ambient$, and real numbers
$\gamma_1,\ldots, \gamma_k$, are given by
%
%
\begin{equation}
\label{XLcharfn}
\E\Biggl\{\!\exp\Biggl\{\!i\sum_{j=1}^k \gamma_jX_L(\bt_j)\!\Biggr\}\!\Biggr\}
= \exp\Biggl\{\!\int_{S_L}\!
\Biggl(\! \exp\Biggl\{\! i\sum_{j=1}^k \gamma_jx(\bt_j)\!\Biggr\} -1\!
\Biggr)
\lambda_X(dx)\! \Biggr\}.\hspace*{-35pt}
\end{equation}
The second, or ``residual,'' component $Y$ has characteristic functions
\begin{eqnarray*}
&&\E\Biggl\{\exp\Biggl\{i\sum_{j=1}^k \gamma_jY(\bt_j)\Biggr\}
\Biggr\}\\
&&\qquad= \exp\Biggl\{ -Q(\gamma_1,\ldots, \gamma_k)
\\
&&\hspace*{54.6pt}{} + \int_{\bbr^\ambient\setminus S_L}
\Biggl( \exp\Biggl\{ i\sum_{j=1}^k \gamma_jx(\bt_j)\Biggr\} -1
-i\sum_{j=1}^k \gamma_j \lt x(\bt_j)\rt\Biggr) \lambda_X(d{x})\\
&&\hspace*{248pt}{} +
iL_1(\gamma_1,\ldots, \gamma_k)\Biggr\},
\end{eqnarray*}
where we are using the notation of
(\ref{ejontchf}), and
\[
L_1(\gamma_1,\ldots, \gamma_k) = L(\gamma_1,\ldots, \gamma_k)
- \int_{S_L} \sum_{j=1}^k \gamma_j\lt x(\bt_j)\rt\lambda
_X(d{x}).
\]
We shall ultimately show that the limiting behaviour of the critical
points of
$X$ depends only on the component $ X^L$, so we study it first.

\item[(iii)] \textit{A limit theorem for the critical points of $X^L$.}
We start by noting that it follows from the form of the characteristic function
(\ref{XLcharfn}) and the definition (\ref{elevyX})
that $X^L$ can, in law, be written as
%
%
\begin{equation}\label{eCompPoissrep}
X^L(\bt) = \sum_{m=1}^N X_mf(S_k, \bt),
\end{equation}
where $N$ is a Poisson random variable with mean $\theta$ given by
(\ref{efinitelevy}), independent of an i.i.d. sequence of random
pairs $((X_m,S_m)), m=1,2,\ldots)$ taking values in
$(\bbr\setminus\{ 0\})\times S$ with the common law $\theta^{-1}F$
restricted to the set
\[
\Bigl\{ (s,x)\in(\bbr\setminus\{ 0\})\times S\dvtx
\sup_{\bt\in\ambient} |xf(s;\bt)|>1\Bigr\}.
\]
Recall that $F$ is the L\'evy measure of the infinitely divisible
random measure $M$ in~(\ref{eidfield}).

Since the sum in (\ref{eCompPoissrep}) is a.s. finite, and the
kernel $f$ has bounded $C^2$ sections $f(s;\cdot)$ for all $s\in
S$, it follows that $X^L$ is bounded and $C^2$ on $\ambient$.

We now decompose the compound Poisson term $X^L$ itself into
a sum of two independent pieces, the stochastically larger of which
will be responsible for the limiting behavior of the critical points of
$ X$. For $u>0$
and $1/2<\beta<1$, define the sequence of independent events
\begin{eqnarray*}
A_m(u) &=& \Bigl\{ \max\Bigl[ \sup_{\bt\in\ambient}|X_mf(S_m;\bt
)|,
\sup_{\bt\in\ambient, i=1,\ldots, d}|X_mf_i(S_m;\bt)|, \\
&&\hspace*{105pt}
\sup_{\bt\in\ambient, i,j=1,\ldots, d}|X_mf_{ij}(S_m;\bt
)|\Bigr] >
u^\beta
\Bigr\}
\end{eqnarray*}
and write
%
%
\begin{eqnarray}
\label{edecomposeXL}
X^L(\bt) &=& \sum_{m=1}^N X_mf(S_m; \bt) \indic_{A_m(u)}
+ \sum_{m=1}^N X_mf(S_m; \bt) \indic_{A_m(u)^c}
\nonumber\\[-8pt]\\[-8pt]
&\definedas& X^{(L,1)}(\bt)+
X^{(L,2)}(\bt). \nonumber
\end{eqnarray}
In Lemma~\ref{ldecomposeXL}, we shall show that $X^{(L,2)}$
and its partial derivatives have
suprema the tail probabilities of which decay faster than the function
$H$, and so are unlikely to affect the critical points of $X$. We shall
return to this point later.\vadjust{\goodbreak}

Now, however, we shall concentrate on the
critical points over high
levels of $X^{(L,1)}$. Define two new events
%
%
\begin{equation}
\label{thebees}\quad
B_1(u) = \Biggl\{ \sum_{m=1}^N \indic(A_m(u)) =1\Biggr\},\qquad
B_2(u) = \Biggl\{ \sum_{m=1}^N \indic(A_m(u)) \geq2\Biggr\}.
\end{equation}
The first of these occurs when there is a single large term in the
Poisson sum (\ref{eCompPoissrep}), the second when there are more.
On the event $B_1(u)$, we define the random
variable $K(u)$ to be the index of large term, and otherwise
allow it to be arbitrarily.

In the notation of Section~\ref{seclimitEuler} in general and Theorem
\ref{tmain} in particular, it follows that, on the event
$B_1(u)$, the following representation holds for the numbers of the critical
points of $X^{(L,1)}$
over the level $u$. For $k=0,1,\ldots, d$, a face $J\in{\mathcal
J}_k$ and
$i=0,1,\ldots, k$,
\begin{eqnarray*}
&&N_{X^{(L,1)}}(J;i:u)
\\[-2pt]
&&\qquad=
\indic\bigl( X_{K(u)}>0\bigr)
\sum_{l=0}^{c_i(J;S_{K(u)})} \indic\bigl(
X_{K(u)}f\bigl( S_{K(u)}; \bt_l\bigl(J;i;S_{K(u)}\bigr)\bigr)>u\bigr)
\\[-2pt]
&&\qquad\quad{} + \indic\bigl( X_{K(u)}<0\bigr)
\sum_{l=0}^{c_{k-i}(J;S_{K(u)})} \indic\bigl(
X_{K(u)}f\bigl( S_{K(u)}; \bt_l\bigl(J;k-i;S_{K(u)}\bigr)\bigr)>u\bigr).
\end{eqnarray*}
Therefore, for any number $r=1,2,\ldots\,$, on the event $B_1(u)$, we
have\break
$ N_{X^{(L,1)}}(J;i:u)\geq r$ if, and only if,
\[
X_{K(u)}>\bigl( f_{[r]}^{(J;i:+)}\bigl(S_{K(u)}\bigr)\bigr)^{-1} u
\quad\mbox{or}\quad X_{K(u)}<-\bigl( f_{[r]}^{(J;k-i:-)}\bigl(S_{K(u)}\bigr)\bigr)^{-1}
u.
\]
We conclude that for any numbers $n(J;i)=1,2,\ldots\,$, for all
$J\in{\mathcal J}_k$, and for all $k=0,1,\ldots, d$ and
$i=0,1,\ldots, k$,
%
%
\begin{eqnarray}\label{eextremeL1}
&&\P\bigl\{ \bigl\{ N_{X^{(L,1)}}(J;i:u)\geq n(J;i)
\mbox{ for all
$J$ and $i$} \bigr\} \cap B_1(u)\bigr\}
\nonumber\\
&&\qquad
= \P\Bigl\{ \Bigl[ X_{K(u)}> \max_{J\in{\mathcal
J}_k, k=0,1,\ldots, d, i=0,1,\ldots, k}\bigl(
f_{[n(J;i)]}^{(J;i:+)}\bigl(S_{K(u)}\bigr)\bigr)^{-1} u
\mbox{ or}
\nonumber\\[-8pt]\\[-8pt]
&&\hspace*{51pt}X_{K(u)}<- \max_{J\in{\mathcal
J}_k, k=0,1,\ldots, d, i=0,1,\ldots, k}\bigl(
f_{[n(J;i)]}^{(J;k-i:-)}\bigl(S_{K(u)}\bigr)\bigr)^{-1} u\Bigr] \nonumber\\
&&\hspace*{273.5pt}{}\cap B_1(u)
\Bigr\}.\nonumber
\end{eqnarray}
Write ${\cal E}_m$ for the union of sets $(-\infty,-\max]$ and
$[\max,\infty)$, where the ``max'' come from the preceding lines with
$K(u)$ replaced by $m$. Then
%
%
\begin{eqnarray}
\label{esplitB1}
&&\P\bigl\{ \bigl\{ X_{K(u)}\in{\cal E}_{K(u)}\bigr\} \cap
B_1(u)\bigr\}
\nonumber\\
&&\qquad= \P\Biggl\{ \bigcup_{m=1}^N \biggl( A_m(u) \cap
\bigcap_{m_1\not= m}
A_{m_1}(u)^c \cap\{ X_m\in{\cal E}_m\}\biggr)\Biggr\}
\nonumber\\[-8pt]\\[-8pt]
&&\qquad= e^{-\theta}\sum_{n=0}^\infty\frac{\theta^n}{n!}
n\P\Biggl\{
{\cal E}_1(u) \cap\bigcap_{m_1=2}^n
A_{m_1}(u)^c \cap\{ X_1\in{\cal E}_1\}\Biggr\}
\nonumber\\
&&\qquad= \theta\P\bigl\{
A_1(u) \cap\{ X_1\in{\cal E}_1\}\bigr\}
- \P\bigl\{ \bigl\{ X_{K(u)}\in{\cal E}_{K(u)}\bigr\} \cap
B_2(u)\bigr\}.\nonumber
\end{eqnarray}
Applying this to the right-hand side of
(\ref{eextremeL1}) and using part (iii) of Lemma~\ref{ldecomposeXL} yields
%
%
\begin{eqnarray}\label{eextremeL11}
&& \P\bigl\{ \bigl\{ N_{X^{(L,1)}}(J;i:u)\geq n(J;i)\mbox{ for all
$J$ and $i$} \bigr\} \cap B_1(u)\bigr\}
\nonumber\\
&&\qquad
= \theta\P\Bigl\{ A_1(u) \cap\Bigl\{ X_1 > \max_{J\in
{\mathcal
J}_k, k=0,1,\ldots, d, i=0,1,\ldots, k}\bigl(
f_{[n(J;i)]}^{(J;i:+)}(S_1)\bigr)^{-1} u\Bigr\} \Bigr\}
\nonumber\\[-8pt]\\[-8pt]
&&\qquad\quad{} + \theta\P\Bigl\{ A_1(u) \cap\Bigl\{ X_1 <- \max
_{J\in{\mathcal
J}_k, k=0,1,\ldots, d, i=0,1,\ldots, k}\bigl(
f_{[n(J;i)]}^{(J;k-i:-)}(S_1)\bigr)^{-1} u\Bigr\} \Bigr\}\hspace*{-4pt}
\nonumber\\
&&\qquad\quad{} - Q_{\mathrm{small}}(u),\nonumber
\end{eqnarray}
where $Q_{\mathrm{small}}(u)/H(u)\to0$ as $u\to\infty$.

Assume for the moment that all the $n(J;i)$ are strictly positive.
Since the parameter $\beta$ in the definition of the event $A_1(u)$ is
less than 1, it follows that, as $u\to\infty$,
\begin{eqnarray*}
&&\P\Bigl\{ A_1(u) \cap\Bigl\{ X_1 > \max_{J\in{\mathcal
J}_k, k=0,1,\ldots, d, i=0,1,\ldots, k}\bigl(
f_{[n(J;i)]}^{(J;i:+)}(S_1)\bigr)^{-1} u\Bigr\} \Bigr\}
\\
&&\qquad
\sim\P\Bigl\{ X_1 > \max_{J\in{\mathcal
J}_k, k=0,1,\ldots, d, i=0,1,\ldots, k}\bigl(
f_{[n(J;i)]}^{(J;i:+)}(S_1)\bigr)^{-1} u\Bigr\}
\\
&&\qquad
= \frac1\theta\int_S \rho\Bigl( s; \Bigl(
\max_{J\in{\mathcal
J}_k, k=0,1,\ldots, d, i=0,1,\ldots, k}\bigl(
f_{[n(J;i)]}^{(J;i:+)}(s)\bigr)^{-1} u, \infty\Bigr)\Bigr)
m(ds).
\end{eqnarray*}
In the last step, we used the law of $X_1$ introduced after
(\ref{eCompPoissrep}) and the decomposition
(\ref{elevymdecompose}) of the measure $F$, and in the middle one
the asymptotic equivalence means that the two ratio of the two
probabilities tends to 1 as $u\to\infty$. Since a similar
asymptotic expression can be written for the second term in the
right-hand side of~(\ref{eextremeL11}), we obtain
\begin{eqnarray*}
&&\lim_{u\to\infty} \frac{\P\{ \{
N_{X^{(L,1)}}(J;i:u)\geq
n(J;i) \mbox{ for all $J$ and $i$} \} \cap B_1(u)\}}{H(u)}
\\
&&\qquad
= \lim_{u\to\infty} H(u)^{-1}\\
&&\qquad\quad{}\times
\int_S \Bigl[ \rho\Bigl( s; \Bigl(
\max_{J\in{\mathcal
J}_k, k=0,1,\ldots, d, i=0,1,\ldots, k}\bigl(
f_{[n(J;i)]}^{(J;i:+)}(s)\bigr)^{-1} u, \infty\Bigr)\Bigr)
\\
&&\hspace*{26pt}\qquad\quad{}
+ \rho\Bigl( s; \Bigl( -\infty,
- \max_{J\in{\mathcal
J}_k, k=0,1,\ldots, d, i=0,1,\ldots, k}\bigl(
f_{[n(J;i)]}^{(J;k-i:-)}(s)\bigr)^{-1} u \Bigr)\Bigr)\Bigr]\\
&&\hspace*{55pt}{}\times m(ds),
\end{eqnarray*}
provided the last limit exists. Applying (\ref{esmallmaxf})
and Potter's bounds, as in Lem\-ma~\ref{lreglevy}, to justify
an interchange of limit and integration, and noting
Assumption~\ref{assregvar} relating $\rho$, $\omega$ and $H$, we have
\begin{eqnarray*}
&&
\lim_{u\to\infty} \frac{\P\{ \{
N_{X^{(L,1)}}(J;i:u)\geq
n(J;i) \mbox{ for all $J$ and $i$} \} \cap
B_1(u)\}}{H(u)}
\\
&&\qquad
= \int_S \Bigl[ w_+(s)\min_{J\in{\mathcal
J}_k, k=0,1,\ldots, d, i=0,1,\ldots, k}\bigl(
f_{[n(J;i)]}^{(J;i:+)}(s)\bigr)^{\alpha}
\\
&&\hspace*{13pt}\qquad\quad{}
+ w_-(s)\min_{J\in{\mathcal
J}_k, k=0,1,\ldots, d, i=0,1,\ldots, k}\bigl(
f_{[n(J;i)]}^{(J;k-i:-)}(s)\bigr)^{\alpha}\Bigr] m(ds)
\\
&&\qquad\definedas
I_{\mathrm c}.
\end{eqnarray*}
Finally, since by part (iii) of Lemma~\ref{ldecomposeXL}, the
event $B_2(u)$ has a probability of a smaller order, we can also
conclude that
%
%
\begin{equation}\label{efinalL1}
\lim_{u\to\infty} \frac{ \P\{ N_{X^{(L,1)}}(J;i:u)\geq
n(J;i) \mbox{ for all $J$ and $i$}\}}{H(u)} = I_{\mathrm c}.
\end{equation}
In view of (\ref{ereglevy}), we can rewrite this as
%
%
\begin{eqnarray}\label{efinalL12}
&&\lim_{u\to\infty} \P\Bigl\{ N_{X^{(L,1)}}(J;i:u)\geq
n(J;i) \mbox{ for all $J$ and $i$} \big|
\sup_{t\in M} X_t\geq u\Bigr\}
\nonumber\\[-8pt]\\[-8pt]
&&\qquad= \frac{I_{\mathrm c}}{\int_S [
w_+(s)\sup_{\bt\in M} f(s,\bt)_+^{\alpha} +
w_-(s)\sup_{\bt\in M} f(s,\bt)_-^{\alpha}] m(ds)}.
\nonumber
\end{eqnarray}

This will complete the proof of the theorem,
at least for the case of strictly positive $n(J;i)$,
once we show that the lighter-tailed
random fields $Y$ of (\ref{erfdecompose}) and $ X^{(L,2)}$
of (\ref{edecomposeXL}) do not change the asymptotic distribution of
the numbers of critical points of $X$.
This will take us a while to show, and makes up the remainder of the proof.

Before we do this, note that handling situations in which
some or all of the numbers $n(J;i)$ are zero is actually only an issue of
semantics, once we recall our convention regarding the 0th order statistic
introduced prior to the statement of the theorem. For example, in the
case when \textit{all} the $n(J;i)$ are zero, the event
on the left-hand side of (\ref{efinalL12}) should be interpreted as
stating that $ X^{(L,1)}$ has crossed the level $u$, given that it has
done so. Not surprisingly, the resulting limit, and the right-hand side,
turn out to be 1. Similar reductions work when only some of the
$n(J;i)$ are zero.

\item[(iv)] \textit{An outline of what remains to do.}
It follows from what we have done so far that
%
%
\begin{equation} \label{esplitX}
X(\bt) = X^{(L,1)}(\bt) +X^{(L,2)}(\bt) +Y(\bt),\qquad \bt\in
\ambient,
\end{equation}
or, equivalently, that
%
%
\begin{equation} \label{esplitXmain}
X^{(L,1)}(\bt) = X(\bt) -X^{(L,2)}(\bt) -Y(\bt),\qquad \bt\in
\ambient.
\end{equation}
What we plan to show is that when when either
$ X$ or $X^{(L,1)}$ reaches a high level~$u$, then the
lighter-tailed
random fields $ Y$ and $X^{(L,2)}$ can be thought of as small perturbations,
both in terms of their absolute values, and those of their first and second
order partial derivatives. This will imply that the asymptotic
conditional joint
distributions of the number of the critical points of the random fields
$X$ and $X^{(L,1)}$ are not affected by the lighter
tailed fields and, hence, coincide.

In fact, what we establish is that near every critical point of one of the
random fields $X$ and $X^{(L,1)}$ there is a critical point, of the same
index, of the other. Equation (\ref{esplitX})
allows us to do this in one direction, and
(\ref{esplitXmain}) will give us the other direction. The
two equations are of the same type, and the fact that the terms in the
right-hand side of (\ref{esplitX}) are independent, while the terms
in the
right-hand side of (\ref{esplitXmain}) are not, will play no role in the
argument. Therefore, we shall treat in detail only one of the two
directions, and describe only briefly the additional steps needed for
the other. The first steps in this program involve collecting some probabilistic
bounds on the closeness of critical points and the behavior of Hessians
there.

\item[(v)] \textit{Bounds on critical points and Hessians.}
We start by introducing a function $D\dvtx S\to(0,\infty]$ that
describes what we
think of as the degree of nondegeneracy of the critical points of an
$s$-section of the kernel $f$. This includes the minimal Euclidian
distance between two distinct critical points of an $s$-section of
the kernel $f$ and the smallest absolute value of an eigenvalue of the
Hessian matrices of the section evaluated at critical points.
Specifically, starting with critical points, and
recalling the definition of the $t_l(J;i;s)$ as the critical
points of index $i$ on the face $J$ for the $s$-section of $X$, define
\begin{eqnarray*}
&&D_1(s) =\min\{ \| \bt_{l_1}(J_1;i_1;s) -
\bt_{l_2}(J_2;i_2;s)\|\dvtx J_j\in{\mathcal J}_{k_j},
\\
&&\hspace*{63pt} 0\leq k_1,k_2\leq d, 0\leq i_j\leq k_j,
0\leq l_j\leq c_{i_j}(J_j;s), j=1,2\},
\end{eqnarray*}
where the minimum is taken over distinct points. Furthermore, define
\begin{eqnarray*}
D_2(s) &=& \min\{ |\lambda|\dvtx\lambda\mbox{ is an
eigenvalue of }
( f_{mn}( s; \bt_{l}(J;i;s)))_{m,n\in\sigma(J)};
\\
&&\hspace*{60.5pt}J\in{\mathcal J}_{k}, 0\leq k \leq d,
0\leq i\leq k,
1\leq l\leq c_{i}(J;s)\}.
\end{eqnarray*}
As usual, both minima are defined to be equal to $+\infty$ if taken
over an empty set.

Now set
%
%
\begin{equation}\label{eDm}
D(s) = \min( D_1(s), D_2(s)).
\end{equation}
Note that, by Assumption~\ref{assregkernel},
$D$ is a strictly positive function,
so that for any any $S$-valued random variable $W$ one has
$\lim_{\tau\to0} \P\{D(W)\leq\tau\} = 0$. Choose $W$ to have the
law $N_W$ given by
%
%
\begin{equation} \label{ecnorm}
\frac{dN_W}{dm}(s) = c_*\bigl( w_+(s)+ w_-(s)\bigr) \sup_{\bt\in
I_d} |f(s,\bt)|^{\alpha},\qquad s\in S,
\end{equation}
where $c_*$ is a normalising constant. That this is possible is a consequence
of (\ref{esmallmaxf}). For $\vep>0$,
choose $\tau_0>0$ so small that $\P\{D(W)\leq\tau_0\}\leq
\vep$. With the random variable $K(u)$ as before, Lemma
\ref{lseparate} gives us that
%
%
\begin{equation} \label{ereduce1}
\limsup_{u\to\infty}\frac{\P\{ \{ D(
S_{K(u)})
\leq\tau_0,
{\sup_{\bt\in I_d}} |X^{(L,1)}(\bt)|>u\} \cap
B_1(u)\}}{H(u)}
\leq c_*^{-1}\vep,\hspace*{-18pt}
\end{equation}
where $B_1(u)$ was defined at (\ref{thebees}) and indicates that there was
only one ``large'' component in the decomposition of $X$.

Note that, since the event
$\{ {\sup_{\bt\in I_d}} |X^{(L,1)}(\bt)|>u\} \cap
B_1(u)$
is a subset of $B_1(u)$, on this event
$X^{(L,1)}(\bt) = X_{K(u)}
f(S_{K(u)}; \bt)$ for all $\bt\in\ambient$. Thus, again on this
this event,
since the supremum of this field
over $I_d$ exceeds $u$, while the kernel $f$ is uniformly
bounded, we conclude that $|X_{K(u)}|>u/\|f\|_\infty$. Therefore, on
the event
\[
\Bigl\{ D\bigl( S_{K(u)}\bigr)
> \tau_0,
\sup_{\bt\in I_d} \bigl|X^{(L,1)}(\bt)\bigr|>u\Bigr\} \cap
B_1(u)
\]
the smallest eigenvalue length
%
%
\begin{eqnarray}
\label{Dmin}
&&D_{\min} \definedas\min\bigl\{ |\lambda|\dvtx\lambda\mbox{
is an
eigenvalue of }
\bigl( X^{(L,1)}_{mn}( \bt)
\bigr)_{m,n\in\sigma(J)};
\nonumber\\[-8pt]\\[-8pt]
&&\hspace*{69pt}
J\in{\mathcal J}_{k}, 0\leq k \leq d,
\mbox{$\bt$ is a critical point on
$J$}\bigr\} \nonumber
\end{eqnarray}
satisfies $D_{\min}> (\tau_0/\|f\|_\infty)u$.

We now combine (\ref{ereduce1}) with (\ref{efinalL1})
as follows. Introduce the event
$\widetilde\Omega_\tau(u)$ that occurs whenever the minimal Euclidian
distance between
two distinct critical points of the random field $(
X^{(L,1)}(\bt), \bt\in I_d)$ is at least $\tau>0$, while
the smallest eigenvalue length of the Hessian evaluated at the critical
points satisfies $D_{\min}> (\tau/\|f\|_\infty)u$. Thus, we have
%
%
\begin{eqnarray}
\label{eL1sep0}\qquad
\liminf_{u\to\infty} \frac{\P\{ \{
N_{X^{(L,1)}}(J;i:u)\geq
n(J;i) \ \forall J,i\} \cap
\widetilde\Omega_{\tau_0}(u)\}}{H(u)}
\geq I_{\mathrm c} - c_*^{-1}\vep,
\end{eqnarray}
where $I_{\mathrm c}$ is as in (\ref{efinalL1}). We can, furthermore,
``sacrifice'' another $\vep$ in the right-hand side of (\ref{eL1sep0})
to add to the event $\widetilde\Omega_\tau(u)$ a requirement that the
\textit{largest} eigenvalue of the Hessian evaluated at the critical
points, which we denote by $D_{\max}$,
satisfies $D_{\max}\leq Mu$ for some positive
$M=M(\epsilon)$. This is possible because $D_{\max}$ is bounded
from above by the largest absolute value of the elements of the Hessian,
which we bound from above by $Mu$ with a large enough $M$. For
the same reason, we can also bound from above the largest value of $\|
\nabla X^{(L,1)}(\bt)\|$ over $I_d$ by $Mu$.

Denoting the resulting event by $\Omega_\tau(u)$, we obtain
%
%
\begin{eqnarray}
\label{eL1sep}\qquad
\liminf_{u\to\infty} \frac{\P\{\{
N_{X^{(L,1)}}(J;i:u)\geq
n(J;i) \ \forall J,i \} \cap
\Omega_{\tau_0}(u)\}}{H(u)}
\geq I_{\mathrm c} - 2c_*^{-1}\vep.
\end{eqnarray}

Now note that since, as stated above, $X^L$ is bounded and $C^2$ on
$\ambient$, and the same is true for $X$ by Assumption~\ref{assc2},
it follows that the\vadjust{\goodbreak} ``remainder'' $Y$ in
(\ref{erfdecompose}) is also a.s. bounded and $C^2$.
Furthermore, by construction, $Y$ and its first and second order
partial derivatives
have L\'evy measures that are supported on
uniformly bounded functions. Consequently, the tail of their absolute
suprema decays exponentially fast; see~\cite{deacosta1980}. In
particular, for $i,j=1,\ldots, d$,
\begin{eqnarray*}
\lim_{u\to\infty} \frac{\P\{ \sup_{\bt\in\ambient}
|Y(\bt)|>u\}}{H(u)} &=& \lim_{u\to\infty} \frac{\P\{
\sup_{\bt\in\ambient}
|Y_i(\bt)|>u\}}{H(u)}
\\
&=& \lim_{u\to\infty} \frac{\P\{ \sup_{\bt\in\ambient}
|Y_{ij}(\bt)|>u\}}{H(u)}
\\ &=& 0.
\end{eqnarray*}
It follows from this, part (ii) of Lemma
\ref{ldecomposeXL}, and the regular variation of $H$, that there
is a function $l(u)\uparrow\infty$ such that $l(u)/u\to0$ as
$u\to\infty$ and
%
%
\begin{eqnarray}
\label{eL2sep}
\lim_{u\to\infty} \frac{\P\{ \sup_{\bt\in\ambient}
|Y(\bt)|>l(u)\}}{H(u)}
&=& \lim_{u\to\infty} \frac{\P\{ \sup_{\bt\in\ambient}
|Y_{ij}(\bt)|>l(u)\}}{H(u)}
\nonumber\\
&=&
\lim_{u\to\infty} \frac{\P\{ \sup_{\bt\in\ambient}
|X^{L,2}(\bt)|>l(u)\}}{H(u)}
\nonumber\\[-8pt]\\[-8pt]
&=& \lim_{u\to\infty} \frac{\P\{ \sup_{\bt\in
\ambient}
|X^{L,2}_{ij}(\bt)|>l(u)\}}{H(u)}\nonumber\\
&=& 0\nonumber
\end{eqnarray}
for $i,j=1,\ldots, d$.

We now combine (\ref{eL1sep}) and
(\ref{eL2sep}) in the following way. Let $ \Omega^{(1)}_{\tau}(u)$ be
the intersection of the event $ \Omega_{\tau}(u)$ with the
complements of all 4 events whose probabilities are displayed in
(\ref{eL2sep}) and set
\[
\Omega_{\mathit{cr}}(u) = \bigl\{
N_{X^{(L,1)}}\bigl(J;i:(1+\tau_2)u\bigr)\geq
n(J;i) \ \forall J,i \bigr\} \cap
\Omega^{(1)}_{\tau_1}(u).
\]
Then, given $0<\vep_1<1$, and using the regular variation of
$H$, we can find $\tau_1,\tau_2>0$ such that
%
%
\begin{equation}
\label{eL3sep}
\liminf_{u\to\infty} \frac{\P\{\Omega_{\mathit{cr}}(u)\}}{H(u)}
\geq(1-\vep_1)I_{\mathrm c}.
\end{equation}

\item[(vi)] \textit{The} (\textit{almost}) \textit{end of the proof.}
Continuing with the above notation, we now
claim that, on the event $\Omega_{\mathit{cr}}(u)$, for $u$ large
enough so that
%
%
\begin{equation}\label{elargeu}
\frac{u}{l(u)} \geq\max\biggl( \frac{8k\| f\|_\infty}{\tau_1},
\frac{4}{\tau_2}\biggr),
\end{equation}
we also have
%
%
\begin{equation} \label{ealsowhole}
N_{X}(J;i:u)\geq n(J;i),\qquad J\in{\mathcal J}_k,
k=0,1,\ldots, d, i=0,1,\ldots, k.
\end{equation}
Note that, once this is established, we shall have
\[
\liminf_{u\to\infty} \frac{\P\{ N_{X}(J;i:u)\geq n(J;i),
J\in
{\mathcal J}, 0\leq i \leq\dim J \}}{H(u)}
\geq(1-\vep_1)I_{\mathrm c}
\]
and, since this holds for all $0<\vep_1<1$, we also have
%
%
\begin{equation} \label{elowerbound}
\liminf_{u\to\infty} \frac{\P\{ N_{X}(J;i:u)\geq n(J;i),
J\in
{\mathcal J}, 0\leq i \leq\dim J \}}{H(u)} \geq I_{\mathrm c}.
\end{equation}

Combining this with (\ref{ereglevy}) gives Theorem~\ref{tmain},
albeit with an inequality rather than an equality in (\ref{elimdist}).

To obtain the opposite inequality
assume that, to the contrary, for some numbers
$n(J;i)$, 
%
%
\begin{equation}\label{etoomuch}
\lim_{n\to\infty} \frac{\P\{ N_{X}(J;i:u_n)\geq n(J;i),
J\in
{\mathcal J}, 0\leq i\leq\dim J \}}{H(u_n)}
> I_{\mathrm c}
\end{equation}
along some sequence $u_n\uparrow\infty$.

Now proceed by repeating the steps performed above and,
this time using (\ref{esplitXmain}) rather than (\ref{esplitX}), and
so demonstrate the existence of a critical point of $X^{(L,1)}$ near
each one of $ X$. Thus, (\ref{etoomuch}) also holds with $X$
replaced by
$X^{(L,1)})$, viz.
\[
\lim_{n\to\infty} \frac{\P\{ N_{X^{(L,1)}}(J;i:u_n)\geq n(J;i),
J\in
{\mathcal J}, 0\leq i\leq\dim J \}}{H(u_n)}
> I_{\mathrm c}.
\]
Since this contradicts (\ref{efinalL1}), (\ref{etoomuch}) cannot be
true, we have the required lower bound, and the proof of Theorem
\ref{tmain} is complete, modulo the need to establish the
claim~(\ref{ealsowhole}).

\item[(vii)] \textit{Establishing} (\ref{ealsowhole}) \textit{to finish the proof.}
In order to establish (\ref{ealsowhole}), we shall
show that, on the event $\Omega_{\mathit{cr}}(u)$, to every
critical point above the level $(1+\tau_2)u$ of the random field $X^{(L,1)}$
we can associate a critical
point above the level $u$ of $X$ which is in the same face
and of the same type.

To this end, let $\bt_0$ be a
critical point above the level $(1+\tau_2)u$ of $X^{(L,1)}$
that belongs to a face $J\in
{\mathcal J}_k$ for some $0\leq k\leq d$, and which is of the type
$i$ for some $0\leq i\leq k$. Let $(\be_1,\ldots, \be_k)$ be an
orthonormal basis
of $\bbr^k$ consisting of normalised eigenvectors of the Hessian matrix
%
%
\begin{equation}
\label{hereweshouldusemanifolds1}
{\mathcal H}^{(L,1)}(\bt_0) = \bigl( X^{(L,1)}_{mn}( \bt
_0)
\bigr)_{m,n\in\sigma(J)},
\end{equation}
and let $\lambda_1, \ldots, \lambda_k$ be
the corresponding eigenvalues. Note that, by the definition of the
event $\Omega_{\mathit{cr}}(u)$, we have
$|\lambda_n|>(\tau_1/\|f\|_\infty)u$ for $n=1,\ldots, k$.
We naturally embed the vectors
$(\be_1,\ldots, \be_k)$ into the face $J$ and make them $d$-dimensional
vectors by appending to them the $d-k$ fixed coordinates of the face
$J$. [We shall continue to denote these vectors by $(\be_1,\ldots,
\be_k)$.]
Note
that for small real numbers $\epsilon_1,\ldots, \epsilon_k$ we have
%
%
\begin{equation}
\label{hereweshouldusemanifolds2}
\nabla X^{(L,1)}\Biggl( \bt_0 + \sum_{j=1}^k \epsilon_j\be_j\Biggr)
= \sum_{j=1}^k \epsilon_j \lambda_j\be_j +
o(\max(|\epsilon_1|,\ldots, |\epsilon_k|)).
\end{equation}
In particular, the directional derivatives
\[
g^{(L,1)}_j(\bt) \definedas\bigl\langle\nabla X^{(L,1)}(\bt),
\be_j\bigr\rangle,\qquad
j=1,\ldots, k,
\]
satisfy
%
%
\begin{equation}
\label{eL4sep}
g^{(L,1)}_j\Biggl( \bt_0 + \sum_{j=1}^k \epsilon_j\be_j\Biggr)
= \epsilon_j\lambda_j + o(\max(|\epsilon_1|,\ldots,
|\epsilon_k|)).
\end{equation}

In what follows, we shall work with a small positive number
$\epsilon>0$, placing more and more conditions on it as we progress, to
clarify precisely how small it will need to be. As a first step,
take $\epsilon<\tau_1/2$, where
$\tau_1$ is as in (\ref{eL3sep}).

Consider a
$k$-dimensional cube (a subset of the face $J$) defined
by
\[
C_\epsilon= \Biggl\{ \bt_0 + \sum_{j=1}^k \theta_j\be_j,
|\theta_j|\leq\epsilon, j=1,\ldots, k\Biggr\},
\]
along with its ($k-1$)-dimensional faces
\[
F_n^{\pm} = \Biggl\{ \bt_0 + \sum_{j=1}^k \theta_j\be_j,
\theta_n = \pm\epsilon, |\theta_j|\leq\epsilon, 1\leq j\leq k,
j\not= n\Biggr\},
\]
where $n=1,\ldots, k$. It follows from (\ref{eL4sep}) that, for
$\epsilon>0$ small enough, $u>1$, and, as above, $M$ large enough, we have
%
%
\begin{equation} \label{eL5sep}
2M\epsilon u \geq2\epsilon|\lambda_n| \geq\bigl| g^{(L,1)}_n(
\bt)\bigr| \geq\frac{\epsilon|\lambda_n|}{2}
\geq\frac{\tau_1 \epsilon}{2\|f\|_\infty}u
\end{equation}
for all $\bt\in F_n^{\pm}$, $n=1,\ldots, k$. The assumption that
$\epsilon$ be small enough now entails that (\ref{eL5sep}) holds
for all
critical points and for all relevant $n$. Since the number of
critical points is finite, this requirement is easy to satisfy.
Similarly, the continuity of the eigenvalues
of a quadratic matrix in its components (see, e.g., Section
7.2. and Corollary 2 in Section 7.4 of~\cite{lancaster1969})
shows that, for all $\epsilon>0$ small enough, the
eigenvalues of the matrix of the second order partial derivatives
$( X_{mn}(\bt)_{m,n\in\sigma(J)}$ have all absolute values
satisfying $|\lambda_n|>(\tau/2\|f\|_\infty)u$ for $n=1,\ldots, k$ and
\mbox{$\bt\in C_\epsilon$}. Finally, we require that $\epsilon$ be small
enough that this lower bounds holds for all
critical points $\bt_0$ considered above. In particular, this
implies that the signs of these
eigenvalues throughout $C_\epsilon$ are the same as those at the point
$\bt_0$.

Next, for a nonempty $I\subset\{ 1,\ldots, k\}$ and $p\in\{-1,1\}^k$
consider the vector
%
%
\begin{equation}\label{evt}
{x}(I,p) = \sum_{i\in I} p_i\be_i.
\end{equation}
Consider a point $\bt$ that belongs to the (relative to the face $J$)
boundary of the cube $C_\epsilon$ and, more specifically, belongs
to the face defined by
%
%
\begin{equation} \label{efaceC}
\biggl( \bigcap_{i\in I, p_i=1}F_i^+\biggr)\cap\biggl( \bigcap_{i\in
I, p_i=-1}F_i^-\biggr)\vadjust{\goodbreak}
\end{equation}
and to no other ($k-1$)-dimensional face of $C_\epsilon$. Define a
function $h^{(L,1)}\dvtx C_\epsilon\to\bbr$ by
\[
h^{(L,1)}(\bt) = \sum_{i=1}^k \bigl( g^{(L,1)}_i(\bt)\bigr)^2.
\]
This is a $C^1$-function, and its gradient (within the face $J$) is
given by
\[
\nabla h^{(L,1)}(\bt) = 2\sum_{i=1}^k g^{(L,1)}_i(\bt) \nabla
g^{(L,1)}_i(\bt)
= 2\sum_{i=1}^k g^{(L,1)}_i(\bt) {\mathcal H}^{(L,1)}(\bt)\be_i^T.
\]
Note also that for all $I$ and $p$ as above,
\[
\bigl\langle\nabla g^{(L,1)}_i(\bt_0), {x}(I,p)\bigr\rangle
= \cases{
\lambda_i p_i, &\quad if $i\in I$,\cr
0, &\quad if $i\notin I$.}
\]
In particular, for any $\bt$ belonging to the face of
$C_\epsilon$ defined by (\ref{efaceC}),
%
%
\begin{eqnarray}
\label{eatt}
&&\bigl\langle\nabla h^{(L,1)}(\bt), {x}(I,p)\bigr\rangle
\nonumber\\
&&\qquad= 2\sum_{i=1}^k g^{(L,1)}_i(\bt) \bigl\langle\nabla
g^{(L,1)}_i(\bt), {x}(I,p)\bigr\rangle
\nonumber\\[-8pt]\\[-8pt]
&&\qquad
= 2\sum_{i\in I} \lambda_i p_i g^{(L,1)}_i(\bt)
\nonumber\\
&&\qquad\quad{} + 2\sum_{i=1}^k g^{(L,1)}_i(\bt)
\bigl\langle\bigl( \nabla
g^{(L,1)}_i(\bt) - \nabla g^{(L,1)}_i(\bt_0) \bigr), {x}(I,p)\bigr\rangle.
\nonumber
\end{eqnarray}
It follows from (\ref{eL4sep}) and (\ref{eL5sep}) that
%
%
\begin{eqnarray}
\label{eL6sep}
&&g^{(L,1)}_n( \bt)>0 \qquad\mbox{for $\bt\in F_n^{+}$ if $\lambda_n>0$
and for $\bt\in F_n^{-}$ if $\lambda_n<0$},
\nonumber\\[-8pt]\\[-8pt]
&&g^{(L,1)}_n( \bt)<0 \qquad\mbox{for $\bt\in F_n^{+}$ if
$\lambda_n<0$
and for $\bt\in F_n^{-}$ if $\lambda_n>0$}.\nonumber
\end{eqnarray}
Consequently, we can conclude, by (\ref{eL6sep}) and (\ref{eL5sep}),
that the first term in the right-hand side of (\ref{eatt}) is
negative and, more specifically, does not exceed
\[
-2\mbox{Card}(I)D_{\min} \frac{\tau_1 \epsilon}{2\|f\|_\infty}u
\leq-(\tau_1/\|f\|_\infty)^2\epsilon u^2.
\]
We can bound the absolute value of the second term in the
right-hand side of (\ref{eatt}) from above by
\[
2k \sum_{i=1}^k \bigl|g^{(L,1)}_i(\bt)\bigr| \cdot\bigl\| \nabla
g^{(L,1)}_i(\bt) - \nabla g^{(L,1)}_i(\bt_0) \bigr\|
\leq2k^2M^2 \epsilon u^2
\]
by the definition of the event $\Omega_{\mathit{cr}}(u)$. This, obviously,
indicates that, for $\epsilon>0$ small enough,
\[
\bigl\langle\nabla h^{(L,1)}(\bt), {x}(I,p)\bigr\rangle\leq-C\epsilon u^2,
\]
where $C$ is a finite positive constant determined by the parameters
in the event $\Omega_{\mathit{cr}}(u)$. Writing
$g_j(\bt) = \langle\nabla X(\bt), \be_j\rangle$, $j=1,\ldots, k$, if
we define
%
%
\begin{equation} \label{eh}
h(\bt) = \sum_{i=1}^k ( g_i(\bt))^2,
\end{equation}
then, on the event $\Omega_{\mathit{cr}}(u)$,
\[
\langle\nabla h(\bt), {x}(I,p)\rangle\leq-C\epsilon u^2 +kl(u)^2.
\]
Taking into account that $l(u)/u\to0$ as $u\to\infty$, where $l$ is
given by (\ref{eL2sep}), we see that
for $u$ large enough it is possible to choose $\epsilon>0$ small
enough such that
%
%
\begin{equation}\label{eatt1}
\langle\nabla h(\bt), {x}(I,p)\rangle< 0
\end{equation}
for any $\bt$ belonging to the face of
$C_\epsilon$ defined by (\ref{efaceC}). The final
requirement on $\epsilon$ is that (\ref{eatt1}) holds.

Similarly, since by the definition of the event $\Omega_{\mathit{cr}}(u)$,
the first order partial derivatives of $X^{(L,2)}$ and $Y$ are bounded by
$l(u)=o(u)$ in
absolute value, we have that (\ref{eL6sep}) and (\ref{elargeu})
also give us
%
%
\begin{eqnarray}
\label{esigng}
&&g_n( \bt)>0 \qquad\mbox{for $\bt\in F_n^{+}$ if $\lambda_n>0$
and for $\bt\in F_n^{-}$ if $\lambda_n<0$},
\nonumber\\[-8pt]\\[-8pt]
&&
g_n( \bt)<0 \qquad\mbox{for $\bt\in F_n^{+}$ if $\lambda_n<0$
and for $\bt\in F_n^{-}$ if $\lambda_n>0$}.\nonumber
\end{eqnarray}

In order to complete the proof and establish (\ref{ealsowhole}) it
suffices to prove that, on $\Omega_{\mathit{cr}}(u)$, $X$ has a
critical point in the cube $C_\epsilon$. If such a critical point
exists, Lem\-ma~\ref{lcritpts} below implies that it will be above the
level $u$ and of the same type as $\bt_0$. Furthermore, these critical
points of $X$ will all be distinct.

To establish the existence of this critical point, note that,
by the continuity of $\nabla X$ and the compactness of $C_\epsilon$, there
is a point $\bt_1$ in $C_\epsilon$ at which the norm of the vector function
$\bg(\bt) = ( g_1(\bt),\ldots, g_k(\bt))$ achieves its minimum
over $C_\epsilon$. We shall prove that, in fact, $\bg(\bt_1)= 0$. By the
linear independence of the basis vectors $\be_1, \ldots, \be_k$, this
will imply that $g_j(\bt_1)=0$ for $j=1,\ldots, k$, and so $\bt_1$ is,
indeed, a~critical point.

Suppose that, to the contrary, $\bg(\bt_1)\not= 0$, and consider firstly
the possibility that the point $\bt_1$ belongs to the (relative to the
face $J$) interior of $C_\epsilon$. Note that the Jacobian of the
transformation $g\dvtx C_\epsilon\to\bbr^k$ is given by
$J_g(\bt) = {\mathcal E}{\mathcal H}(\bt)$,
where ${\mathcal H}(\bt) = ( X_{mn}( \bt)
)_{m,n\in\sigma(J)}$ is the Hessian of
$X$, and ${\mathcal E}$ is a $k\times k$ matrix with rows $\be_1,
\ldots, \be_k$. We have already established above
that, on the event $\Omega_{\mathit{cr}}(u)$,
${\mathcal H}$ is nondegenerate throughout $C_\epsilon$. Since the
vectors $\be_1, \ldots, \be_k$ are linearly independent, we conclude
that the matrix ${\mathcal E}$ is nondegenerate as well. Since the
vector $\bg(\bt_1)$ does not vanish, it has a nonvanishing component.
Without loss\vadjust{\goodbreak} of generality, we can assume that $ g_1(\bt_1)\not=0$.
Choose a
vector ${x}\in\bbr^k$ for which
$J_g(\bt_1) {x}' = (1,0,\ldots,0)'.
$
Then for $\delta\in\bbr$, with $|\delta|$ small,
\[
\bg(\bt_1 + \delta{x})
= \bg(\bt_1) + \delta J_g(\bt_1) {x}^T + o(|\delta|)
\]
and so
\[
\| \bg(\bt_1 + \delta{x})\|^2
=\sum_{j=1}^k g_j(\bt_1)^2 + 2\delta g_1(\bt_1)+ o(|\delta|)
< \sum_{j=1}^k g_j(\bt_1)^2 = \| \bg(\bt_1
)\|^2
\]
for $\delta$ with $|\delta|$ small enough and such that $\delta
g_1(\bt_1)<0$. This contradicts the assumed minimality of $\|
\bg(\bt_1)\|$ and so we must have $g(t_1)=0$, as
required, for this case.

It remains to consider the case $\bg(\bt_1)\not= 0$, but the point
$\bt_1$ belongs to the boundary of the cube $C_\epsilon$. Let
$\bg(\bt_1)$ belong to the face of the cube defined by
(\ref{efaceC}). With the function $h$ defined in (\ref{eh}), we
have, for $\delta>0$ small,
\[
h\bigl( \bt_1 + \delta{x}(I,p)\bigr) =
h( \bt_1) + \delta( \nabla h(\bt), {x}(I,p))
+ o(\delta).
\]
By (\ref{eatt1}), this last expression is smaller than $h(
\bt_1)$ if $\delta>0$ is small enough. However, by the
definition of the vector ${x}(I,p)$, the point $ \bt_1 + \delta
{x}(I,p)$ belongs to $C_\epsilon$ for $\delta>0$ small. Once again,
this contradicts the assumed minimality of $\|
\bg(\bt_1)\|$.

Thus, we have established (\ref{ealsowhole}) and,
therefore, (\ref{elowerbound}), and so the theorem, modulo the
need to prove the following three lemmas. 
\end{longlist}
%
%
\begin{Lemma} \label{ldecomposeXL} The following three results hold:

\begin{longlist}[(iii)]
\item[(i)]
The random fields $X^{(L,1)}$ and $X^{(L,2)}$
on the right-hand side of the decomposition (\ref{edecomposeXL}) are
independent.

\item[(ii)] The random field
$ X^{(L,2)}$
has $C^2$ sample functions and satisfies
\begin{eqnarray*}
\lim_{u\to\infty} \frac{\P\{ \sup_{\bt\in\ambient}
|X^{(L,2)}(\bt)|>u\}}{H(u)}
&=& \lim_{u\to\infty} \frac{\P\{ \sup_{\bt\in\ambient}
|X_i^{(L,2)}(\bt)|>u\}}{H(u)}
\\
&=& \lim_{u\to\infty} \frac{\P\{ \sup_{\bt\in\ambient}
|X_{ij}^{(L,2)}(\bt)|>u\}}{H(u)}
=0,
\end{eqnarray*}
$i,j=1,\ldots, d$,
where $H$ is the regularly varying function of Assumption
\ref{assregvar}.

\item[(iii)]
The number of terms in the sum defining $X^{(L,1)}$ satisfies
\[
\lim_{u\to\infty}\frac{\P\{ \sum_{m=1}^N \indic(A_m(u))
\geq2\}}
{H(u)} = 0.
\]
\end{longlist}
\end{Lemma}
\begin{pf}
The claim (i) follows from the fact that a
Poisson random measure, when restricted to disjoint measurable sets,
forms independent Poisson random measures on these sets (see,
e.g.,~\cite{resnick1992}). Since the sum defining the random field
$ X^{(L,2)}$ is a.s. finite, the fact
that it has sample functions in $C^2$ follows from Assumption\vadjust{\goodbreak}
\ref{assc2}. Furthermore, for $\epsilon>0$, choose
$n_\epsilon>0$ so large that $\P\{N>n_\epsilon\}\leq\epsilon$. The
above discussion implies that
the number $K(u,\epsilon)$ of the terms in the sum defining
$X^{(L,2)}$ in (\ref{edecomposeXL}) that satisfy
\[
\sup_{\bt\in\ambient}|X_mf(S_m;\bt)|>\frac{u}{2n_\epsilon}
\]
is Poisson with the mean no greater than
\begin{eqnarray*}
&& F\Bigl\{ (s,x)\in(\bbr\setminus\{ 0\})\times S\dvtx
\sup_{\bt\in\ambient} |xf(s;\bt)|> u/(2n_\epsilon) \Bigr\}
\\
&&\qquad
= \lambda_X\Bigl\{ g\dvtx\sup_{\bt\in\ambient}
|g(\bt)|>u/(2n_\epsilon)\Bigr\} \sim CH(u)
\end{eqnarray*}
as $u\to\infty$, where we have used (\ref{elevyX}) and Lemma
\ref{lreglevy}. Therefore, for large $u$
\[
\P\Bigl\{ \sup_{\bt\in\ambient} \bigl|X^{(L,2)}(\bt)\bigr|>u \Bigr\}
\leq\P\{N>n_\epsilon\} + \P\{K(u,\epsilon)\geq2\}
\leq\epsilon+ CH(u)^2
\]
and so
\[
\limsup_{u\to\infty} \frac{\P\{ \sup_{\bt\in\ambient}
|X^{(L,2)}(\bt)|>u\}}{H(u)}\leq\epsilon.
\]
Letting $\epsilon\to0$ completes the proof of the first limit in part
(ii) of the lemma, and the other limits are established in the same
way. Part (iii) of the lemma can be proven similarly.
\end{pf}
%
%
\begin{Lemma} \label{lseparate} With the $c_*$ of (\ref{ecnorm}),
the $D_m$ in (\ref{eDm}) satisfy
\[
\limsup_{u\to\infty}\frac{\P\{ \{ D(
S_{K(u)})
\leq\tau_0,
\sup_{\bt\in I_d} |X^{(L,1)}(\bt)|>u\} \cap
B_1(u)\}}{H(u)}
\leq c_*^{-1}\vep.
\]
\end{Lemma}
\begin{pf}
We use a decomposition as in (\ref{esplitB1}) to obtain
\begin{eqnarray*}
&&\P\Bigl\{ \Bigl\{ D\bigl( S_{K(u)}\bigr)\leq\tau_0, \sup
_{\bt\in
I_d} \bigl|X^{(L,1)}(\bt)\bigr|>u\Bigr\} \cap B_1(u)\Bigr\}
\\
&&\qquad\leq\theta\P\Biggl\{ A_1(u) \cap\bigcap_{m_1=2}^n
A_{m_1}(u)^c \cap\Bigl\{ D(S_1)\leq\tau_0,
\sup_{\bt\in I_d} \bigl|X^{(L,1)}(\bt)\bigr|>u\Bigr\}\Biggr\}.
\end{eqnarray*}
Since on the event $A_{m}(u)^c$ one has $\sup_{\bt\in
\ambient}|X_mf(S_m;\bt)|\leq u^\beta$, it follows that the latter
probability can be asymptotically bounded by
\begin{eqnarray*}
&&\theta\P\Bigl\{ D(S_1)\leq\tau_0, \sup_{\bt\in
I_d}|X_1f(S_1;\bt)|>u \Bigr\}
\\
&&\qquad= \int_S\int_{\bbr\setminus\{ 0\}} \indic\Bigl( D(s)\leq
\tau_0,
|x|\sup_{\bt\in I_d}|f(s;\bt)|>u\Bigr) F(ds,dx)
\\ &&\qquad= \int_S \indic\bigl( D(s)\leq\tau_0\bigr) \biggl( \int
_{\bbr\setminus\{
0\}} \indic\Bigl( |x|> u\Bigl(\sup_{\bt\in
I_d}|f(s;\bt)|\Bigr)^{-1}\Bigr)
\rho(s; dx)\biggr) m(ds)
\\ &&\qquad\sim H(u) \int_S \indic\bigl( D(s)\leq\tau_0\bigr)
\bigl( w_+(s)+ w_-(s)\bigr) \sup_{\bt\in
I_d}|f(s;\bt)|^\alpha m(ds),
\end{eqnarray*}
where we used Assumptions~\ref{assregvar},~\ref{assregkernel} and
(\ref{esmallmaxf}). The
lemma now follows from the choice of $\tau_0$.
\end{pf}
%
%
\begin{Lemma} \label{lcritpts}
Suppose that for every critical point $\bt_0$ of the random field
$(
X^{L,1}(\bt), \bt\in I_d )$,
the random field $(
X(\bt), \bt\in I_d )$ has, on the event $\Omega_{\mathit{cr}}(u)$, a
critical point in the cube $C_\epsilon$. Then the
critical points of $X$ in $ I_d$ correspond to distinct critical
points of $X^{L,1}$, are themselves distinct,
are all above the level $u$, and each of them is of the same type as the
corresponding critical point of $X^{L,1}$.
\end{Lemma}
\begin{pf}
The fact that the
critical points of $X$ corresponding to distinct critical
points of the field $X^{L,1}$ are all distinct follows from
the lower bound on the distance between two distinct critical points
of $X^{L,1}$ in the definition of the
event $\Omega_{\mathit{cr}}(u)$ and the choice of $\epsilon$. The
fact that
all the critical points are above the level $u$ follows from the lower
bounds on the values of $X^{L,1}$ at its critical points in the
definition of $\Omega_{\mathit{cr}}(u)$ and, once again, the choice of
$\epsilon$. It remains, therefore, to prove that a critical point in
the cube $C_\epsilon$ of $X$ is of the same type as the
critical point $\bt_0$ of $X^{L,1}$.

To this end, note that the absolute values of the
eigenvalues of the matrix of the second order partial derivatives
$( X^{L,1}_{mn}(\bt) + Y_{mn}(\bt))_{m,n\in\sigma(J)}$ are,
on the event $\Omega_{\mathit{cr}}(u)$,
bounded from above by $2k l(u)$. Using continuity of the eigenvalues
of a quadratic matrix in its components (see, once again, Section
7.2 and Corollary~2 in Section 7.4 of~\cite{lancaster1969}), we see
that the Euclidian distance between an eigenvalue
of
$( X_{mn}(\bt))_{m,n\in\sigma(J)}$ and the corresponding
eigenvalue of
$( X^{L,1}_{mn}(\bt))_{m,n\in\sigma(J)}$ is bounded from
above by $2k l(u)$. Using the choice of $\epsilon$ then shows that the
numbers of the negative eigenvalues of the two Hessians
are identical, as required.
\end{pf}

\printaddresses


\begin{thebibliography}{20}

\bibitem{adlertaylor2007}
\begin{bbook}[mr]
\bauthor{\bsnm{Adler},~\bfnm{Robert~J.}\binits{R.~J.}} \AND
  \bauthor{\bsnm{Taylor},~\bfnm{Jonathan~E.}\binits{J.~E.}}
(\byear{2007}).
\btitle{Random Fields and Geometry}.
\bpublisher{Springer}, \baddress{New York}.
\bid{mr={2319516}}
\bptok{imsref}%
\end{bbook}
\endbibitem

\bibitem{adlertaylorworsley2010}
\begin{bbook}[author]
\bauthor{\bsnm{Adler},~\bfnm{R.~J.}\binits{R.~J.}},
  \bauthor{\bsnm{Taylor},~\bfnm{J.~E.}\binits{J.~E.}} \AND
  \bauthor{\bsnm{Worsley},~\bfnm{K.~J.}\binits{K.~J.}}
(\byear{2013}).
\btitle{Applications of Random Fields and Geometry: Foundations and Case
  Studies}.
\bpublisher{Springer}, \baddress{New York}.
\bnote{Early chapters available at
\href{http://webee.technion.ac.il/people/adler/publications.html}{webee.technion.ac.il/people/adler/publications.html}}.
\bptok{imsref}%
\end{bbook}
\endbibitem

\bibitem{binghamgoldieteugels1987}
\begin{bbook}[mr]
\bauthor{\bsnm{Bingham},~\bfnm{N.~H.}\binits{N.~H.}},
  \bauthor{\bsnm{Goldie},~\bfnm{C.~M.}\binits{C.~M.}} \AND
  \bauthor{\bsnm{Teugels},~\bfnm{J.~L.}\binits{J.~L.}}
(\byear{1987}).
\btitle{Regular Variation}.
\bseries{Encyclopedia of Mathematics and Its Applications}
\bvolume{27}.
\bpublisher{Cambridge Univ. Press}, \baddress{Cambridge}.
\bid{mr={0898871}}
\bptok{imsref}%
\end{bbook}
\endbibitem

\bibitem{bravermansamorodnitsky1998a}
\begin{barticle}[mr]
\bauthor{\bsnm{Braverman},~\bfnm{Michael}\binits{M.}} \AND
  \bauthor{\bsnm{Samorodnitsky},~\bfnm{Gennady}\binits{G.}}
(\byear{1998}).
\btitle{Symmetric infinitely divisible processes with sample paths in {O}rlicz
  spaces and absolute continuity of infinitely divisible processes}.
\bjournal{Stochastic Process. Appl.}
\bvolume{78}
\bpages{1--26}.
\bid{doi={10.1016/S0304-4149(98)00055-6}, issn={0304-4149}, mr={1653284}}
\bptok{imsref}%
\end{barticle}
\endbibitem

\bibitem{deacosta1980}
\begin{barticle}[mr]
\bauthor{\bparticle{de} \bsnm{Acosta},~\bfnm{Alejandro}\binits{A.}}
(\byear{1980}).
\btitle{Exponential moments of vector-valued random series and triangular
  arrays}.
\bjournal{Ann. Probab.}
\bvolume{8}
\bpages{381--389}.
\bid{issn={0091-1798}, mr={0566603}}
\bptok{imsref}%
\end{barticle}
\endbibitem

\bibitem{lancaster1969}
\begin{bbook}[mr]
\bauthor{\bsnm{Lancaster},~\bfnm{Peter}\binits{P.}}
(\byear{1969}).
\btitle{Theory of Matrices}.
\bpublisher{Academic Press}, \baddress{New York}.
\bid{mr={0245579}}
\bptok{imsref}%
\end{bbook}
\endbibitem

\bibitem{linde1986}
\begin{bbook}[mr]
\bauthor{\bsnm{Linde},~\bfnm{Werner}\binits{W.}}
(\byear{1986}).
\btitle{Probability in {B}anach Spaces---Stable and Infinitely Divisible
  Distributions},
\bedition{2nd} ed.
\bpublisher{Wiley}, \baddress{Chichester}.
\bid{mr={0874529}}
\bptok{imsref}%
\end{bbook}
\endbibitem

\bibitem{maejimarosinski2002}
\begin{barticle}[mr]
\bauthor{\bsnm{Maejima},~\bfnm{Makoto}\binits{M.}} \AND
  \bauthor{\bsnm{Rosi{\'n}ski},~\bfnm{Jan}\binits{J.}}
(\byear{2002}).
\btitle{Type {$G$} distributions on {$\Bbb R\sp d$}}.
\bjournal{J. Theoret. Probab.}
\bvolume{15}
\bpages{323--341}.
\bid{doi={10.1023/A:1015044726122}, issn={0894-9840}, mr={1898812}}
\bptok{imsref}%
\end{barticle}
\endbibitem

\bibitem{marcusrosinski2005}
\begin{barticle}[mr]
\bauthor{\bsnm{Marcus},~\bfnm{Michael~B.}\binits{M.~B.}} \AND
  \bauthor{\bsnm{Rosi{\'n}ski},~\bfnm{Jan}\binits{J.}}
(\byear{2005}).
\btitle{Continuity and boundedness of infinitely divisible processes: A
  {P}oisson point process approach}.
\bjournal{J. Theoret. Probab.}
\bvolume{18}
\bpages{109--160}.
\bid{doi={10.1007/s10959-004-2579-1}, issn={0894-9840}, mr={2132274}}
\bptok{imsref}%
\end{barticle}
\endbibitem

\bibitem{Potter}
\begin{barticle}[mr]
\bauthor{\bsnm{Potter},~\bfnm{H.~S.~A.}\binits{H.~S.~A.}}
(\byear{1940}).
\btitle{The mean values of certain {D}irichlet series, {II}}.
\bjournal{Proc. Lond. Math. Soc. (2)}
\bvolume{47}
\bpages{1--19}.
\bid{issn={0024-6115}, mr={0005141}}
\bptok{imsref}%
\end{barticle}
\endbibitem

\bibitem{rajputrosinski1989}
\begin{barticle}[mr]
\bauthor{\bsnm{Rajput},~\bfnm{Balram~S.}\binits{B.~S.}} \AND
  \bauthor{\bsnm{Rosi{\'n}ski},~\bfnm{Jan}\binits{J.}}
(\byear{1989}).
\btitle{Spectral representations of infinitely divisible processes}.
\bjournal{Probab. Theory Related Fields}
\bvolume{82}
\bpages{451--487}.
\bid{doi={10.1007/BF00339998}, issn={0178-8051}, mr={1001524}}
\bptok{imsref}%
\end{barticle}
\endbibitem

\bibitem{resnick1992}
\begin{bbook}[mr]
\bauthor{\bsnm{Resnick},~\bfnm{Sidney}\binits{S.}}
(\byear{1992}).
\btitle{Adventures in Stochastic Processes}.
\bpublisher{Birkh\"auser}, \baddress{Boston, MA}.
\bid{mr={1181423}}
\bptok{imsref}%
\end{bbook}
\endbibitem

\bibitem{rosinski1989}
\begin{barticle}[mr]
\bauthor{\bsnm{Rosi{\'n}ski},~\bfnm{Jan}\binits{J.}}
(\byear{1989}).
\btitle{On path properties of certain infinitely divisible processes}.
\bjournal{Stochastic Process. Appl.}
\bvolume{33}
\bpages{73--87}.
\bid{doi={10.1016/0304-4149(89)90067-7}, issn={0304-4149}, mr={1027109}}
\bptok{imsref}%
\end{barticle}
\endbibitem

\bibitem{rosinski1990b}
\begin{barticle}[mr]
\bauthor{\bsnm{Rosi{\'n}ski},~\bfnm{Jan}\binits{J.}}
(\byear{1990}).
\btitle{On series representations of infinitely divisible random vectors}.
\bjournal{Ann. Probab.}
\bvolume{18}
\bpages{405--430}.
\bid{issn={0091-1798}, mr={1043955}}
\bptok{imsref}%
\end{barticle}
\endbibitem

\bibitem{rosinski1991}
\begin{bincollection}[mr]
\bauthor{\bsnm{Rosi{\'n}ski},~\bfnm{Jan}\binits{J.}}
(\byear{1991}).
\btitle{On a class of infinitely divisible processes represented as mixtures of
  {G}aussian processes}.
In \bbooktitle{Stable Processes and Related Topics ({I}thaca, {NY}, 1990)}
(\beditor{\bfnm{S.}\binits{S.}~\bsnm{Cambanis}},
\beditor{\bfnm{G.}\binits{G.}~\bsnm{Samorodnitsky}} \AND
\beditor{\bfnm{M.~S.}\binits{M.~S.}~\bsnm{Taqqu}}, eds.).
\bseries{Progress in Probability}
\bvolume{25}
\bpages{27--41}.
\bpublisher{Birkh\"auser}, \baddress{Boston, MA}.
\bid{mr={1119350}}
\bptok{imsref}%
\end{bincollection}
\endbibitem

\bibitem{rosinskisamorodnitsky1993}
\begin{barticle}[mr]
\bauthor{\bsnm{Rosi{\'n}ski},~\bfnm{Jan}\binits{J.}} \AND
  \bauthor{\bsnm{Samorodnitsky},~\bfnm{Gennady}\binits{G.}}
(\byear{1993}).
\btitle{Distributions of subadditive functionals of sample paths of infinitely
  divisible processes}.
\bjournal{Ann. Probab.}
\bvolume{21}
\bpages{996--1014}.
\bid{issn={0091-1798}, mr={1217577}}
\bptok{imsref}%
\end{barticle}
\endbibitem


\bibitem{samorodnitskytaqqu1994}
\begin{bbook}[mr]
\bauthor{\bsnm{Samorodnitsky},~\bfnm{Gennady}\binits{G.}} \AND
  \bauthor{\bsnm{Taqqu},~\bfnm{Murad~S.}\binits{M.~S.}}
(\byear{1994}).
\btitle{Stable Non-{G}aussian Random Processes: Stochastic Models With Infinite Variance}.
\bpublisher{Chapman \& Hall}, \baddress{New York}.
\bid{mr={1280932}}
\bptok{imsref}%
\end{bbook}
\endbibitem

\bibitem{sato1999}
\begin{bbook}[mr]
\bauthor{\bsnm{Sato},~\bfnm{Ken-iti}\binits{K.-i.}}
(\byear{1999}).
\btitle{L\'evy Processes and Infinitely Divisible Distributions}.
\bseries{Cambridge Studies in Advanced Mathematics}
\bvolume{68}.
\bpublisher{Cambridge Univ. Press}, \baddress{Cambridge}.
\bid{mr={1739520}}
\bptok{imsref}%
\end{bbook}
\endbibitem

\end{thebibliography}
\end{document}